\documentclass[aos,seceqn,citesort,dvips]{arximspdf}
\usepackage{mathbh}
\usepackage{graphicx}

\doi{10.1214/10-AOS800}
\volume{38}
\issue{4}
\pubyear{2010}
\firstpage{2422}
\lastpage{2464}

\makeatletter

\newproclaim{ass}{Assumption}
\newtheorem{theo}{Theorem}
\newtheorem{prop}{Proposition}
\newtheorem{lemma}{Lemma}

\makeatother

\begin{document}
\begin{frontmatter}

\title{A deconvolution approach to estimation of a~common shape in a
shifted curves model}
\runtitle{Deconvolution for shifted curves}

\begin{aug}
\author[A]{\fnms{J\'{e}r\'{e}mie} \snm{Bigot}\ead[label=e1]{Jeremie.Bigot@math.univ-toulouse.fr}\corref{}} and
\author[A]{\fnms{S\'{e}bastien} \snm{Gadat}\ead[label=e2]{Sebastien.Gadat@math.univ-toulouse.fr}}
\runauthor{J. Bigot and S. Gadat}
\affiliation{Universit\'{e} de Toulouse}
\address[A]{Institut de Math\'{e}matiques de Toulouse et CNRS\\
Universit\'{e} de Toulouse\\
31062 Toulouse Cedex 9\\
France\\
\printead{e1}\\
\phantom{E-mail: }\printead*{e2}} 
\end{aug}

\pdfauthor{Jeremie Bigot, Sebastien Gadat}

\received{\smonth{4} \syear{2009}}
\revised{\smonth{1} \syear{2010}}

%
\begin{abstract}
This paper considers the problem of adaptive estimation of a mean
pattern in a randomly shifted curve model. We show that this problem
can be transformed into a linear inverse problem, where the density of
the random shifts plays the role of a convolution operator. An
adaptive estimator of the mean pattern, based on wavelet thresholding
is proposed. We study its consistency for the quadratic risk as the
number of observed curves tends to infinity, and this estimator is
shown to achieve a near-minimax rate of convergence over a large class
of Besov balls. This rate depends both on the smoothness of the common
shape of the curves and on the decay of the Fourier coefficients of the
density of the random shifts. Hence, this paper makes a connection
between mean pattern estimation and the statistical analysis of linear
inverse problems, which is a new point of view on curve registration
and image warping problems. We also provide a new method to estimate
the unknown random shifts between curves. Some numerical experiments
are given to illustrate the performances of our approach and to compare
them with another algorithm existing in the literature.
\end{abstract}

%

\end{frontmatter}

\section{Introduction}
\subsection{Model and objectives}

In many fields of interests including biology, medical imaging, or
chemistry, observations are coming from $n$ individuals curves or
graylevel images. Such observations are commonly referred to as
functional data, and models involving such data have been recently
extensively studied in statistics (see \cite{ramsil,ramsil02}
for a detailed introduction to functional data analysis). In such
settings, it is reasonable to assume that the data at hand
$Y_{m},m=1,\ldots,n$, satisfy the following white noise regression
model:
%
%
\begin{equation} \label{gen:model}
dY_m(x) = f_{m}(x) \,dx + \varepsilon_{m} \,dW_m(x),\qquad x \in\Omega,
m=1,\ldots,n,
\end{equation}
where $\Omega$ is a subset of ${\mathbb{R}}^{d}$, $f_{m} \dvtx\Omega\to
{\mathbb{R}}$
are unknown regression functions, and $W_{m}$ are independent standard
Brownian motions on $\Omega$ with $\varepsilon_{m}$ representing
different levels of additive noise.
In many situations, the individual curves or images have a certain
common structure which may lead to the assumption that they are
generated from some semi-parametric model of the form
%
%
\begin{equation} \label{eq:modsemipar}
f_{m}(x) = f(x,\tau_{m})\qquad\mbox{for } x \in\Omega\mbox{ and
some } \tau_{m} \in{\mathcal T}\subset{\mathbb{R}}^{p},
\end{equation}
where $f \dvtx\Omega\times{\mathcal T}\to{\mathbb{R}}$ represents an
unknown shape
common to all the $f_{m}$'s. This shape function (also called mean
pattern) may depend on unknown individual \textit{random} parameters
$\tau
_{m},m=1,\ldots,n$, belonging to a compact set ${\mathcal T}$ of
${\mathbb{R}}^{p}$,
which model individual variations. Such a semi-parametric
representation for the $f_{m}$'s is the so-called \textit{self-modeling
regression framework} (SEMOR) introduced by Kneip and Gasser \cite{kg}. \textit{Shape
invariant models} (SIM) are a special class of such models for which
(see, e.g., \cite{kg})
%
%
\begin{equation} \label{eq:SIM}
f_{m}(x) = f(\phi(x,\tau_{m})),
\end{equation}
where for any $ \tau\in{\mathcal T}$, the function $x \mapsto\phi
(x,\tau
)$ is a smooth diffeomorphism of $\Omega$ and $\phi\dvtx\Omega\times
{\mathcal T}\to\Omega$ is a \textit{known} function. Models such as
(\ref
{eq:SIM}) are useful to account for shape variability in time between
curves (see, e.g., \cite{gerga,liumuller}) or in space between
images, which is the well-known problem of curve registration or image
warping (see \cite{glasbey} and the discussion therein for an
overview, \cite{BGL,BGV} and references therein). SIM models
(\ref{eq:SIM}) also represent a large class of statistical models to
study the difficult problem of recovering a mean pattern from a set of
similar curves or images in the presence of random deformations and
additive noise, which corresponds to the general setting of Grenander's
theory of shapes \cite{Gre}. The overall objective of this paper is to
discuss the fundamental problem of estimating of the mean pattern $f$
which can then be used to learn nonlinear modes of variations in time
or shape between similar curves or images.

\subsection{Previous work on mean pattern estimation}

Very few results exist in the literature on nonparametric estimation of
$f$ for SIM models (\ref{eq:SIM}) based on noisy data from (\ref
{gen:model}). The problem of estimating the common shape of a set of
curves that differ only by a time transformation is usually referred to
as the curve registration problem in statistics, and it has received a
lot of attention over the last two decades; see, for example,
\cite{big,gaskneip,kneipgas,liumuller,ramli,wanggas}.
However, in these papers, an asymptotic study
as the number of curves $n$ grows to infinity is generally not
considered. Estimation of the shape function for SEMOR models related
to (\ref{gen:model}) and (\ref{eq:modsemipar}) is studied in
\cite{kg} with a double asymptotic in the number $n$ of curves and the
number of observed time points per curve. In the simplest case of
shifted curves, various approaches have been developed. Based on a
model with a fixed number $n$ of curves, semiparametric estimation of
the deformation parameters $\tau_{m}$ and nonparametric estimation of
the shape function is proposed in \cite{mazaloubgam} and
\cite{vimond}. A generalization of this approach for the estimation of
scaling, rotation and translation parameters for two-dimensional images
is proposed in \cite{BGV}. Estimation of a common shape for randomly
shifted curves and asymptotic in $n$ is also considered in
\cite{ronn}. There is also a huge literature in image analysis on mean
pattern estimation, and some papers have recently addressed the problem
of estimating the common shape of a set of similar images with
asymptotic in the number of images; see, for example,
\cite{AAT,BGL,ma} and references therein. However, in all the
above cited papers rates of convergence and optimality of the proposed
estimators for $f$ have not been studied.

\subsection{A benchmark model for nonparametric estimation of a mean pattern}

The simplest SIM model is the case of randomly shifted curves, namely
\[
f_{m}(x) = f(x - \tau_{m})\qquad \mbox{for } x \in[0,1] \mbox{
and } \tau_{m} \in{\mathbb{R}},
\]
that has recently received some attention in the statistical literature
\cite{loubcast,mazaloubgam,ronn,vimond}. In
this paper, it will thus be assumed that we observe realizations of $n$
noisy and randomly shifted curves $Y_1, \ldots, Y_n$ coming from the
following Gaussian white noise model:
%
%
\begin{equation}\label{model}\quad
dY_m(x) = f(x-\tau_m) \,dx + \varepsilon_{m} \,dW_m(x),\qquad x \in[0,1],
m=1,\ldots,n,
\end{equation}
where $f$ is the unknown mean pattern of the curves, $W_{m}$ are
independent standard Brownian motions on $[0,1]$, the $\varepsilon_{m}$'s
represent levels of noise which may vary from curve to curve, and the
$\tau_m$'s are unknown random shifts independent of the $W_{m}$'s. The
aim of this paper is to study some statistical aspects related to the
problem of estimating $f$, and to propose new methods of estimation.

Model (\ref{model}) is realistic in many situations where it is
reasonable to assume that the observed curves represent replications of
almost the same process and when a large source of variation in the
experiments is due to transformations of the time axis. Such a model is
commonly used in many applied areas dealing with functional data such
as neuroscience \cite{IRT} or biology \cite{ronn}. More generally,
the model (\ref{model}) represents a kind of benchmark model for
studying the problem of recovering the mean pattern $f$ in SIM models.
The results derived in this paper show that the model (\ref{model}),
although simple, already provides some new insights on the statistical
aspects of mean pattern estimation.

The function $f \dvtx{\mathbb{R}}\to{\mathbb{R}}$ is assumed to be
periodic with period
$1$, and the shifts $\tau_m$ are supposed to be independent and
identically distributed (i.i.d.) random variables with density $g \dvtx
{\mathbb{R}}
\to{\mathbb{R}}^{+}$ with respect to the Lebesgue measure $dx$ on
${\mathbb{R}}$. Our
goal is to estimate nonparametrically the shape function $f$ on $[0,1]$
as the number of curves $n$ goes to infinity.

Let $L^{2}([0,1])$ be the space of squared integrable functions on
$[0,1]$ with respect to $dx$, and denote by $\|f\|^{2} = \int_{0}^{1}
|f(x)|^{2}\,dx$ the squared norm of a function $f$. Assume that
${\mathcal F}
\subset L^{2}([0,1])$ represents some smoothness class of functions
(e.g., a~Sobolev or a Besov ball), and let $\hat{f}_{n} \in
L^{2}([0,1])$ be some estimator of the common shape $f$, that is, a
measurable function of the random processes $Y_{m}, m=1,\ldots,n$.
For some $f \in{\mathcal F}$, the risk of the estimator $\hat{f}_{n}$ is
defined to be
\[
{\mathcal R}(\hat{f}_{n},f) = {\mathbb E}\|\hat{f}_{n} - f\|^{2},
\]
where the above expectation ${\mathbb E}$ is taken with respect to the
law of
$\{ Y_{m},m=1,\ldots,n \}$. In this paper, we propose to investigate
the optimality of an estimator by introducing the following minimax risk
\[
{\mathcal R}_{n}({\mathcal F}) = \inf_{\hat{f}_{n}} \sup_{f \in
{\mathcal F}} {\mathcal R}(\hat{f}_{n},f),
\]
where the above infimum is taken over the set of all possible
estimators in model~(\ref{model}). One of the main contributions of
this paper is to derive asymptotic lower and upper bounds for
${\mathcal R}
_{n}({\mathcal F})$ which, to the best of our knowledge, has not been
considered before.

Indeed, we show that there exist constants $M_{1}, M_{2}$, a sequence
of reals $r_{n} = r_{n}({\mathcal F})$ tending to infinity, and an estimator
$\hat{f}_{n}^{\ast}$ such that
\[
\lim_{n \to+ \infty} r_{n} {\mathcal R}_{n}({\mathcal F}) \geq M_{1}
\quad\mbox{and}\quad
\lim_{n \to+ \infty} r_{n} \sup_{f \in{\mathcal F}} {\mathcal
R}(\hat{f}_{n}^{\ast
},f) \leq M_{2}.
\]
However, the construction of $\hat{f}_{n}^{\ast}$ may depend on
unknown quantities such as the smoothness of $f$, and such estimates
are therefore called nonadaptive. Since it is now recognized that
wavelet decomposition is a powerful tool to derive adaptive estimators
(see, e.g., \cite{DJKP95jrssb}),\vspace*{1pt} a second contribution of this
paper is thus to propose wavelet-based estimators $\hat{f}_{n}$ that
attain a near-minimax rate of convergence in the sense there exits a
constant $M_2$ such that
\[
\lim_{n \to+ \infty} (\log n)^{-\beta} r_{n} \sup_{f \in{\mathcal
F}} {\mathcal R}
(\hat{f}_{n},f) \leq M_{2}\qquad \mbox{for some } \beta> 0.
\]

\subsection{Main result}

Minimax risks will be derived under particular smoothness assumptions
on the density $g$. The main result of this paper is that the
difficulty of estimating $f$ is quantified by the decay to zero of the
Fourier coefficients $\gamma_{\ell}$ of the density $g$ of the shifts
defined as
%
%
\begin{equation}\label{eq:defgamma}
\gamma_{\ell} = \mathbb{E} ( e^{- i 2 \pi\ell\tau} )
= \int_{-\infty}^{+\infty} e^{- i 2 \pi\ell x} g(x) \,dx
\end{equation}
for $\ell\in{\mathbb{Z}}$. Depending how fast these Fourier
coefficients tend
to zero as $|\ell| \to+ \infty$, the reconstruction of $f$ will be
more or less accurate. This comes from the fact that the expected value
of each observed process $Y_{m}(x)$ is given by
\[
{\mathbb E}Y_{m}(x) = {\mathbb E}f(x-\tau_{m}) = \int_{-\infty
}^{+\infty}
f(x-\tau) g(\tau) \,d \tau\qquad \mbox{for } x \in[0,1].
\]
This expected value is thus the convolution of $f$ by the density $g$
which makes the problem of estimating $f$ an inverse problem whose
degree of ill-posedness and associated minimax risk depend on the
smoothness assumptions on $g$.

This phenomenon is a well-known fact in deconvolution problems (see,
e.g.,
\cite{JKPR04jrssb,sapatinaspensky,PV99aos}), and more
generally for linear inverse problems as studied in~\cite{cavgopitsy}.
In this paper, the following type of assumption on $g$ is considered.
\begin{ass}\label{assordi}
The Fourier coefficients of $g$ have a polynomial decay, that is,
for some real $\nu> 0$, there exist two constants $C_{\max} \geq
C_{\min} > 0$ such that
$
C_{\min} |\ell|^{-\nu} \leq|\gamma_{\ell}| \leq C_{\max} |\ell
|^{-\nu}
$
for all $\ell\in{\mathbb{Z}}$.
\end{ass}

In standard inverse problems, such as deconvolution, the optimal rate
of convergence we can expect from an arbitrary estimator typically
depends on such smoothness assumptions. The parameter $\nu$ is usually
referred to as the degree of ill-posedness of the inverse problem, and
it quantifies the difficult of inverting the convolution operator. The
following theorem shows that a similar phenomenon holds for the minimax
risk associated to model (\ref{model}). Note that to simplify the
presentation, all the theoretical results are given for the simple
setting where the level of noise is the same for all curves, \textit
{that is},
$\varepsilon_{m} = \varepsilon$ for all $m=1,\ldots,n$ and
some $\varepsilon> 0$. Finally, one also needs the following assumption
on the decay of the density~$g$.
\begin{ass} \label{assdecay}
There exists a constant $C > 0$ and a real $\alpha> 1$ such that the
density $g$ satisfies
$g(x) \leq\frac{C}{1+ |x|^{\alpha}} \mbox{ for all } x \in{\mathbb{R}}$.
\end{ass}

Note that Assumption \ref{assdecay} is not a very restrictive
condition as $g$ is supposed to be an integrable function on ${\mathbb{R}}$.
\begin{theo} \label{theo:main}
Suppose that the smoothness class ${\mathcal F}$ is a Besov ball
$B_{p,q}^{s}(A)$ of radius $A > 0$ with $p,q \geq1$ and smoothness
parameter $s > 0$ (a precise definition of Besov spaces will be given
later on). Suppose that $g$ satisfies Assumptions \ref{assordi}
and~\ref{assdecay}. Let $p'=\min(2,p)$ and assume that $s \geq1/p'$. If
$s > 2\nu+ 1$, then
\[
r_{n}({\mathcal F}) = n^{({2s})/({2s + 2 \nu+1})}.
\]
\end{theo}

Hence, Theorem \ref{theo:main} shows that under Assumption \ref
{assordi} the minimax rate $r_{n}$ is of polynomial order of the sample
size $n$, and that this rate deteriorates as the degree of
ill-posedness $\nu$ increases. Such a behavior is well known for
standard periodic deconvolution in the white noise model
\cite{JKPR04jrssb,sapatinaspensky}, and Theorem \ref{theo:main}
shows that a similar phenomenon holds for the model (\ref{model}). To
the best of our knowledge, this is a new result which makes a
connection between mean pattern estimation and the statistical analysis
of deconvolution problems.

\subsection{Fourier analysis and an inverse problem formulation}

Let us first remark that the model (\ref{model}) exhibits some
similarities with periodic deconvolution in the white noise model as
described in \cite{JKPR04jrssb}. For $x \in[0,1]$, let us define the
following density function:
%
%
\begin{equation} \label{eq:G}
G(x) = \sum_{k \in{\mathbb{Z}}} g(x+k).
\end{equation}
Note that $G(x)$ exists for all $x \in[0,1]$ provided $g$ has a
sufficiently fast decay at infinity; see Assumption \ref{assdecay}.
Since $f$ is periodic with period 1, one has
\[
\int_{-\infty}^{+\infty} f(x-\tau) g(\tau) \,d \tau= \int_{0}^{1}
f(x-\tau) G(\tau) \,d \tau,
\]
and note that $\gamma_{\ell} = \int_{-\infty}^{+\infty} e^{- i 2
\pi\ell x} g(x) \,dx = \int_{0}^{1} e^{- i 2 \pi\ell x} G(x) \,dx $.
Hence, if one defines $\xi_{m}(x) = f(x-\tau_{m}) - \int_{0}^{1}
f(x-\tau) G(\tau) \,d\tau$ and $\xi(x) = \frac{1}{n} \sum_{m=1}^{n}
\xi_{m}(x)$, then taking the mean of the $n$ equations in (\ref
{model}) yields the model
%
%
\begin{eqnarray} \label{model2}
dY(x) &=& \int_{0}^{1} f(x-\tau) G(\tau) \,d \tau \,dx +
\xi(x)\,dx\nonumber\\[-8pt]\\[-8pt]
&&{} + \frac
{\varepsilon}{\sqrt{n}} \,dW(x),\qquad x \in[0,1],\nonumber
\end{eqnarray}
with
%
%
\begin{equation}\label{eq:esp}
\varepsilon^{2} = \frac{1}{n} \sum_{m=1}^{n} \varepsilon_{m}^{2},
\end{equation}
and where $W(x)$ is a standard Brownian motion $[0,1]$.

The model (\ref{model2}) differs from the periodic deconvolution model
investigated in \cite{JKPR04jrssb} by the error term $ \xi(x)$.
Asymptotically, $ \xi(x)$ is a Gaussian variable, so this suggests to
use the wavelet thresholding procedures developed in
\cite{JKPR04jrssb} to derive upper bounds for the minimax risk.
However, it
should be noted that the additive error term $\xi(x)$ significantly
complicates the estimating procedure as \textit{the variance of $\xi(x)$
clearly depends on the unknown function $f$}. Moreover, deriving lower
bounds for the minimax risk in models such as (\ref{model2}) is
significantly more difficult than in the standard white noise model
without the additive term $\xi(x)$.

Now let us formulate models (\ref{model}) and (\ref{model2}) in the
Fourier domain. Supposing that $f \in L^{2}([0,1])$, we denote by
$\theta_\ell$ its Fourier coefficients for $\ell\in{\mathbb{Z}}$, namely
$\theta_{\ell} = \int_{0}^1 e^{- 2 i \ell\pi x} f(x) \,dx$.
The model (\ref{model}) can then be rewritten as
%
%
\begin{eqnarray} \label{cml}
c_{m,\ell} :\!&=& \int_{0}^1 e^{- 2 i \ell\pi x} \,dY_{m}(x)  =  \theta
_{\ell} e^{-i 2 \pi\ell\tau_m} + \varepsilon_{m} z_{\ell,m}
\nonumber\\[-8pt]\\[-8pt]
& = & \theta_{\ell} \gamma_{\ell} + \xi_{\ell,m} + \varepsilon_{m}
z_{\ell,m} \nonumber
\end{eqnarray}
with $\xi_{\ell,m} = \theta_{\ell} e^{-i 2 \pi\ell\tau_m} -
\theta_{\ell} \gamma_{\ell} $, and where $z_{\ell,m}$ are i.i.d.
$N_{{\mathbb{C}}} (0,1 )$ variables, that is, complex Gaussian
variables with zero mean and such that ${\mathbb E}|z_{\ell,n}|^{2} = 1$.

Thus, we can compute the sample mean $\tilde{c}_{\ell}$ of the $\ell$th
Fourier coefficient over the $n$ curves as
%
%
\begin{equation} \label{eqseqmodel}
\tilde{c}_{\ell} = \frac{1}{n} \sum_{m=1}^n c_{\ell,m} = \theta
_{\ell} \tilde{\gamma}_{\ell} + \varepsilon\eta_{\ell} = \theta
_{\ell} \gamma_{\ell} + \xi_{\ell} + \varepsilon\eta_{\ell}
\end{equation}
with $ \tilde{\gamma}_{\ell} = \frac{1}{n} \sum_{m=1}^n e^{-i 2
\pi\ell\tau_m}$, $\xi_{\ell} = \frac{1}{n} \sum_{m=1}^{n} \xi
_{\ell,m}$ and where the $\eta_{\ell}$'s are i.i.d. complex Gaussian
variables with zero mean and such that ${\mathbb E}|\eta_{\ell}|^{2}
= \frac
{1}{n}$. The average Fourier coefficients $\tilde{c}_{\ell}$ in
equation (\ref{eqseqmodel}) can thus be viewed as a set of
observations which is very close to a sequence space formulation of a
statistical inverse problem as described, for example, in
\cite{cavgopitsy}. As in model (\ref{model2}), the additive error term
$\xi
_{\ell}$ is asymptotically Gaussian, however its variance is $\frac
{1}{n} | \theta_{\ell}|^{2}(1-|\gamma_{\ell}|^{2})$ which is
obviously unknown as it depends on $f$.

If we assume that the density $g$ of the random shifts is known, one
can perform a deconvolution step by taking
%
%
\begin{equation}\label{eq:deftheta}
\hat{\theta}_{\ell} = \frac{\tilde{c}_{\ell}}{\gamma_{\ell}} =
\theta_{\ell} \frac{ \tilde{\gamma}_{\ell} }{\gamma_{\ell}} +
\varepsilon\frac{ \eta_{\ell} }{\gamma_{\ell}}
\end{equation}
to estimate the Fourier coefficients of $f$ since, for large $n$,
$\frac{\tilde{\gamma}_{\ell}}{\gamma_{\ell}}$ is close to $1$ by
the strong law of large numbers.

Based on the $\hat{\theta}_{\ell}$'s, two types of estimators are
studied. The simplest one uses spectral cut-off with a cutting
frequency depending on the smoothness assumptions on $f$, and is thus
nonadaptive. The second estimator is based on wavelet thresholding and
is shown to be adaptive using the procedure developed in
\cite{JKPR04jrssb}. Note that part of our results are presented for the case
where the coefficients $\gamma_{\ell}$ are known. Such a framework is
commonly used in nonparametric regression and inverse problems to
obtain consistency results and to study asymptotic rates of
convergence, where it is generally supposed that the law of the
additive error is Gaussian with zero mean and \textit{known} variance
$\varepsilon^{2}$; see, for example,
\cite{JKPR04jrssb,sapatinaspensky,cavgopitsy}. In model (\ref{model}), the
random shifts may be viewed as a second source of noise and for the
theoretical analysis of this problem the law of this other random noise
is also supposed to be known.

\subsection{An inverse problem with unknown operator}

If the density $g$ is unknown, one can view the problem of estimating
$f$ in model (\ref{model}) as a deconvolution problem with unknown
eigenvalues which complicates significantly the estimation procedure.
Such a framework corresponds to the general setting of an inverse
problem with a partially unknown operator. Recently, some papers have
addressed this problem; see, for example,
\cite{cavraim,EfromovichK01,hoffreiss,neu97}, assuming that an
independent sample of noisy eigenvalues or noisy operator is available
which allows an estimation of the $\gamma_{\ell}$'s. However, such an
assumption is not applicable to our model (\ref{model}). Therefore, we
introduce a new method for estimating $f$ is the case of an unknown
density $g$ which leads to a new class of estimators to recover a mean pattern.

\subsection{Organization of the paper}

In Section \ref{sec:estimflinear}, we consider a linear but
nonadaptive estimator based on spectral cut-off. In Section \ref
{sec:estimfnonlinear}, a nonlinear and adaptive estimator based on
wavelet thresholding is studied in the case of known density $g$, and
upper bound for the minimax risk are studied over Besov balls. In
Section \ref{sec:minimax}, we derive lower bounds for the minimax
risk. In Section \ref{sec:estimgamma}, it is explained how one can
estimate the mean pattern $f$ when the density $g$ is unknown. Finally,
in Section \ref{sec:simu}, some numerical examples are proposed to
illustrate the performances of our approach and to compare them with
another algorithm proposed in the literature. All proofs are deferred
to a technical \hyperref[app]{Appendix} at the end of the paper.

\section{Linear estimation of the common shape and upper bounds for
the risk for Sobolev balls} \label{sec:estimflinear}

\subsection{Risk decomposition}

For $\ell\in{\mathbb{Z}}$, a linear estimator of the $\theta_{\ell
}$'s is
given by
$\hat{\theta}_{\ell}^{\lambda} = \lambda_{\ell} \frac{\tilde
{c}_{\ell}}{\gamma_{\ell}}$,
where $\lambda= (\lambda_{\ell})_{\ell\in{\mathbb{Z}}}$ is a\vspace*{-1pt}
sequence of
nonrandom weights called a filter. An estimator $\hat{f}_{n,\lambda}$
of $f$ is then obtained via the inverse Fourier transform
$\hat{f}_{n,\lambda}(x) = \sum_{\ell\in{\mathbb{Z}}} \hat{\theta
}_{\ell
}^{\lambda} e^{- i 2 \pi\ell x}$,
and thanks to the Parseval's relation, the risk of this estimator is
given by
${\mathcal R}(\hat{f}_{n,\lambda},f) = \mathbb{E} \sum_{\ell\in
{\mathbb{Z}}} |\hat
{\theta}_{\ell} - \theta_{\ell}|^{2}$.
The problem is then to choose the sequence $(\lambda_{\ell})_{\ell
\in{\mathbb{Z}}}$ in an optimal way. The following proposition gives the
bias-variance decomposition of ${\mathcal R}(\hat{f}_{n,\lambda},f) $.
\begin{prop} \label{prop:decomp}
For any given nonrandom filter $\lambda$, the risk of the estimator
$\hat{f}_{n,\lambda}$ can be decomposed as
%
%
\begin{equation} \label{eq:decomprisk}\qquad
{\mathcal R}(\hat{f}_{n,\lambda},f) =
{\underbrace{ \sum_{\ell\in\mathbb{Z}} (\lambda_{\ell}-1)^2
|\theta_{\ell}|^2}_{\mathit{Bias}}} + {\underbrace{ \frac{1}{n}
\sum_{\ell\in\mathbb{Z}}
\lambda_{\ell}^2 \biggl[|\theta_{\ell}|^2 \biggl(\frac{1}{|\gamma
_{\ell}|^2} -1 \biggr)
+ \frac{\varepsilon^2}{|\gamma_{\ell}|^2} \biggr]}_{\mathit{Variance}}}.
\end{equation}
\end{prop}

Note that the decomposition (\ref{eq:decomprisk}) does not correspond
exactly to the classical bias-variance decomposition for linear inverse
problems. Indeed, the variance term in (\ref{eq:decomprisk}) differs
from the classical expression of the variance for linear estimator in
statistical inverse problems which would be in our notation $\varepsilon
^2 \sum_{\ell\in{\mathbb{Z}}} \frac{\lambda_{\ell}^2}{|\gamma
_{\ell
}|^2} $. Hence, contrary to classical inverse problems, the variance
term of the risk depends also on the Fourier coefficients $\theta
_{\ell}$ of the unknown function $f$ to recover.

\subsection{Linear estimation}

Let us introduce the following smoothness class of functions which can
be identified with a periodic Sobolev ball:
\[
H_{s}(A) = \biggl\{f \in L^{2}([0,1]) ; \sum_{\ell\in{\mathbb{Z}}}
(1+|\ell|^{2s}) |\theta_{\ell}|^{2} \leq A \biggr\}
\]
for some constant $A$ and some smoothness parameter $s > 0$, where
$\theta_{\ell} = \int_{0}^1 e^{- 2 i \ell\pi x} f(x) \,dx$.
Now consider a linear estimator obtained by spectral cut-off,
that is, for a projection filter of the form $\lambda^{M}_{\ell
} =
\mathbh{1}
_{|\ell| \leq M}$ for some integer $M$. For an appropriate choice of
$M$, the following proposition gives the asymptotic behavior of the
risk ${\mathcal R}(\hat{f}_{n,\lambda^{M}},f)$.
\begin{prop}\label{prop:risklinearpoly}
Assume that $f$ belongs to $H_{s}(A)$ for some real $s > 1/2$ and $A >
0$, and that $g$ satisfies Assumption \ref{assordi}, that is, polynomial
decay of the $\gamma_{\ell}$'s. Then, if $M = M_{n}$ is chosen such
that $M_{n} \sim n^{{1}/({2s+ 2\nu+1})}$, then there exists a
constant $C$ not depending on $n$ such that as $n \to+ \infty$
\[
\sup_{f \in H_{s}(A) }{\mathcal R}(\hat{f}_{n,\lambda^{M}},f) \leq C
n^{-({2s})/({2s + 2 \nu+1})}.
\]
\end{prop}

The above choice for $M_{n}$ depends on the smoothness $s$ of the
function $f$ which is generally unknown in practice and such a spectral
cut-off estimator is thus called nonadaptive. Moreover, the result is
only suited for smooth functions since Sobolev balls $H_{s}(A)$ for $s
> 1/2$ are not suited to model shape functions $f$ which may have
singularities such as points of discontinuity.

\section{Nonlinear estimation with Meyer wavelets and upper bounds for
the risk for Besov balls} \label{sec:estimfnonlinear}

Wavelets have been successfully used for various inverse problems
\cite{donoho}, and for the specific case of deconvolution Meyer
wavelets, a
special class of band-limited functions introduced in \cite{M92}, have
recently received special attention in nonparametric regression: see
\cite{JKPR04jrssb} and \cite{sapatinaspensky}.

\subsection{Wavelet decomposition and the periodized Meyer wavelet basis}

This wavelet basis is derived through the periodization of the Meyer
wavelet basis of $L^{2}({\mathbb{R}})$ (see \cite{JKPR04jrssb} for further
details on its construction). Denote by $\phi_{j,k}$ and $\psi_{j,k}$
the Meyer scaling and wavelet functions at scale $j \geq0$ and
location $0 \leq k \leq2^{j}-1$. For any function $f$ of
$L^{2}([0,1])$, its wavelet decomposition can be written as:
$f = \sum_{k =0}^{2^{j_{0}}-1} c_{j_{0},k} \phi_{j_{0},k} + \sum_{j =
j_{0}}^{+ \infty} \sum_{k=0}^{2^{j}-1} \beta_{j,k} \psi_{j,k}$,
where $c_{j_{0},k} = \int_{0}^{1} f(x) \phi_{j_{0},k}(x) \,dx$, $\beta
_{j,k} = \int_{0}^{1} f(x),\psi_{j,k}(x) \,dx$ and $j_{0} \geq0$
denotes the usual coarse level of resolution. Moreover, the squared
norm of $f$ is given by
$\|f\|^{2} = \sum_{k =0}^{2^{j_{0}}-1} c_{j_{0},k}^{2} + \sum_{j =
j_{0}}^{+ \infty} \sum_{k=0}^{2^{j}-1} \beta_{j,k}^{2}$.
It is well known that Besov spaces can be characterized in terms of
wavelet coefficients (see, e.g., \cite{JKPR04jrssb}). Let $s >
0$ denote the usual smoothness parameter, then for the Meyer wavelet
basis and for a Besov ball $B^{s}_{p,q}(A)$ of radius $A > 0$ with $1
\leq p,q \leq\infty$, one has that
$
B^{s}_{p,q}(A) = \{ f \in L^{2}([0,1]) \dvtx (\sum_{k
=0}^{2^{j_{0}}-1} |c_{j_{0},k}|^{p} )^{{1}/{p}}
+ ( \sum_{j = j_{0}}^{+ \infty} 2^{j(s + 1/2-1/p)q}
( \sum_{k=0}^{2^{j}-1} |\beta_{j,k}|^{p} )^{{q}/{p}}
)^{{1}/{q}} \leq A \}
$
with the respective above sums replaced by maximum if $p=\infty$ or
$q=\infty$. 

Meyer wavelets can be used to efficiently compute the coefficients
$c_{j,k}$ and $\beta_{j,k}$ by using the Fourier transform. Indeed,
thanks to the Plancherel's identity, one obtains that
%
%
\begin{equation} \label{eq:plancherel}
\beta_{j,k}= \sum_{\ell\in\Omega_{j}} \psi^{j,k}_{\ell} \theta
_{\ell},
\end{equation}
where $\psi^{j,k}_{\ell} = \int_{0}^{1} \psi_{j,k}(x) e^{-i 2 \pi
\ell x} \,dx $ denote the Fourier coefficients of $\psi_{j,k}$ and
$\Omega_{j} = \{\ell\in{\mathbb{Z}}; \psi^{j,k}_{\ell} \neq0 \}
$. As
Meyer wavelets $\psi_{j,k}$ are band-limited, $\Omega_{j}$ is a
finite subset set of $[-2^{j+2}c_{0},-2^{j}c_{0}] \cup
[2^{j}c_{0},2^{j+2}c_{0}]$ with $c_{0} = 2\pi/3$ (see
\cite{JKPR04jrssb}), and fast algorithms for computing the above sum have
been proposed in \cite{K94} and \cite{rast07}. The coefficients
$c_{j_{0},k} $ can be computed analogously with $\phi$ instead of
$\psi$ and $\tilde{\Omega}_{j_{0}} = \{\ell\in{\mathbb{Z}}; \phi
^{j_{0},k}_{\ell} \neq0 \}$ instead of $\Omega_{j}$.

Hence, the noisy Fourier coefficients $\hat{\theta}_{\ell}$ given by
(\ref{eq:deftheta}) can be used to quickly compute the following
empirical wavelet coefficients of $f$ as
%
%
\begin{equation} \label{eq:empwavcoef}
\hat{\beta}_{j,k} = \sum_{\ell\in\Omega_{j}} \psi^{j,k}_{\ell}
\hat{\theta}_{\ell} \quad\mbox{and}\quad \hat{c}_{j_{0},k} = \sum_{\ell
\in\Omega_{j_{0}}} \phi^{j_{0},k}_{\ell} \hat{\theta}_{\ell}.
\end{equation}

\subsection{Nonlinear estimation via hard-thresholding}

It is well known that adaptivity can be obtained by using nonlinear
estimators based on appropriate thresholding of the estimated wavelet
coefficients (see, e.g., \cite{DJKP95jrssb}). A nonlinear estimator by
hard-thresholding is defined by
%
%
\begin{equation} \label{eq:hatfnonlin}
\hat{f}_{n}^{h} = \sum_{k =0}^{2^{j_{0}}-1} \hat{c}_{j_{0},k} \phi
_{j_{0},k} + \sum_{j = j_{0}}^{j_{1}} \sum_{k=0}^{2^{j}-1} \hat
{\beta}_{j,k} \mathbh{1}_{\{| \hat{\beta}_{j,k} | \geq\lambda
_{j,k}\}
} \psi_{j,k},
\end{equation}
where the $\lambda_{j,k}$'s are appropriate thresholds (positive
numbers), and $j_{1}$ is the finest resolution level used for the
estimator. As shown by \cite{JKPR04jrssb}, for periodic deconvolution
the choice for $j_{1}$ and the thresholds $\lambda_{j,k}$ typically
depends on the degree $\nu$ of ill-posedness of the problem. Following
Theorem 1 in \cite{JKPR04jrssb}, to derive rate of convergence for
$\hat{f}_{n}^{h}$ one has to control moments of order $2$ and $4$ of
$|\hat{\beta}_{j,k} - \beta_{j,k}|$ and the probability of deviation
of $\hat{\beta}_{j,k}$ from $\beta_{j,k}$.
\begin{prop} \label{prop:moments} Let $1 \leq p \leq\infty$, $1 \leq
q \leq\infty$, $s - 1/p +1/2 > 0$ and $A > 0$. Assume that $g$
satisfies Assumptions \ref{assordi} and \ref{assdecay}. Then there
exist positive constants $C_{3}$ and $C_{4}$ such that for any $j \geq
j_{0} \geq0$, $0 \leq k \leq2^{j}-1$ and all $f \in B^{s}_{p,q}(A)$,
${\mathbb E}|\hat{c}_{j_{0},k} - c_{j_{0},k} |^{2} \leq C_{3} \frac
{2^{2j_{0} \nu}}{n} , {\mathbb E}|\hat{\beta}_{j,k} - \beta_{j,k}
|^{2} \leq C_{3} \frac{2^{2j \nu}}{n}$
and
${\mathbb E}|\hat{\beta}_{j,k} - \beta_{j,k} |^{4} \leq C_{4} (
\frac
{2^{j 4 \nu}}{n^{2} } + \frac{2^{j (4 \nu+ 1)}}{n^{3} } )$.
\end{prop}

Proposition \ref{prop:moments} shows that the variance of the
empirical wavelet coefficients is proportional to $ \frac{2^{2j \nu
}}{n}$ which comes from the amplification of the noise by the inversion
of the convolution operator. The choice of the threshold $\lambda
_{j,k}$ is done by controlling the probability of deviation of the
empirical wavelet coefficients $\hat{\beta}_{j,k}$ from the true
wavelet coefficient $\beta_{j,k}$ which is given by the following proposition.
\begin{prop} \label{prop:assproba} Let $f \in B_{p,q}^{s}(A)$, $n
\geq1$ and $j \geq0$. Suppose that $g$ satisfies Assumption \ref
{assdecay}. Let $\eta> 0$. For $j \geq0$ and $0 \leq k \leq2^{j}-1$
define the following threshold:
%
%
\begin{equation} \label{eq:hatthr}
\lambda_{j,k} = \lambda_{j} = 2 \sigma_{j} \sqrt{\frac{2 \eta\log
(n)}{n}}
\end{equation}
with
$\sigma^{2}_{j} = 2^{-j} \varepsilon^{2} \sum_{\ell\in\Omega_{j}} |
\gamma_{\ell} |^{-2}$.
Then, for all sufficiently large $j$,
\[
{\mathbb P} ( |\hat{\beta}_{j,k} - \beta_{j,k} | \geq\lambda_{j}
) \leq2 n^{-\eta}.
\]
\end{prop}

Note that the level-dependent threshold (\ref{eq:hatthr}) corresponds
to the usual universal thresholds for deconvolution problem based on
wavelet decomposition as in \cite{JKPR04jrssb}. Then, using the
thresholds $\lambda_{j}$, we finally arrive at the following theorem
which gives an upper bound for the minimax risk over a large class of
Besov balls.
\begin{theo} \label{theo:upperbound}
Assume that $g$ satisfies Assumptions \ref{assordi} and \ref
{assdecay}. Let $j_{1}$ and $j_{0}$ be the largest integers such that
$
2^{j_{1}} \leq( n \log(n)^{-1} )^{{1}/({2 \nu+1})}
$
and
$2^{j_{0}} \leq\log( \log(n) )$.
Let $\hat{f}_{n}^{h}$ be the nonlinear estimator obtained by
hard-thresholding with the above choice for $j_{1}$ and $j_{0}$, and
using the thresholds $\lambda_{j}$ defined by (\ref
{eq:hatthr}) with $\eta\geq2$. Let $1 \leq p \leq\infty$, $1 \leq q
\leq\infty$ and $A > 0$. Let $ p'=\min(2,p)$, $s' = s + 1/2-1/p$,
and assume that $s \geq1/p'$.

If $s \geq(2 \nu+1)(1/p-1/2)$, then
\[
\sup_{f \in B_{p,q}^{s}(A)} \| \hat{f}_{n}^{h} - f\|^{2} = \mathcal O
\bigl( n ^{-({2s})/({2s + 2\nu+1}) } (\log n)^{\beta} \bigr)
\qquad\mbox{with } \beta= \frac{2s}{2s + 2\nu+1}.
\]
If $s < (2 \nu+1)(1/p-1/2)$, then
\[
{\sup_{f \in B_{p,q}^{s}(A)}} \| \hat{f}_{n}^{h} - f\|^{2} = \mathcal O
\bigl( ( n^{-1}\log(n) )^{({2s'})/({2s'+2\nu})} \bigr).
\]
\end{theo}

In standard periodic deconvolution in the white noise model (see, e.g.,
\cite{JKPR04jrssb}), there exist two different upper bounds for the
minimax rate which are usually referred to as the dense case [$s \geq
(2 \nu+1)(1/p-1/2)$] when the hardest functions to estimate are spread
uniformly over $[0,1]$, and the sparse case [$s < (2 \nu+1)(1/p-1/2)$]
when the worst functions to estimate have only one nonvanishing
wavelet coefficient. Theorem \ref{theo:upperbound} shows that a
similar phenomenon holds for the model (\ref{model}), and to the best
of our knowledge, this is a new result.
%

\section{Minimax lower bound} \label{sec:minimax}

The following theorem gives an asymptotic lower bound on the minimax
risk ${\mathcal R}_{n}(B^{s}_{p,q}(A))$ for a large class of Besov balls.
\begin{theo} \label{theo:lowerbound}
Let $1 \leq p \leq\infty$, $1 \leq q \leq\infty$ and $A > 0$.
Suppose that $g$ satisfies Assumption \ref{assordi}. Let $p' = \min
(2,p)$. Assume that $s \geq1/p'$ and $\nu> 1/2$.

If $s \geq(2 \nu+1)(1/p-1/2)$ and $s > 2\nu+ 1$ (dense
case), there exits a constant $M_{1}$ depending only on $A,s,p,q$ such that
\[
{\mathcal R}_{n}(B^{s}_{p,q}(A)) \geq M_{1} n^{-({2s})/({2s + 2 \nu+
1})}\qquad
\mbox{as } n \to+ \infty.
\]
\end{theo}

In the dense case, the hardest functions to estimate are spread
uniformly over the interval $[0,1]$, and the proof is based on an
adaptation of Assouad's cube technique (see, e.g., Lemma 10.2 in
\cite{HKPT}) to the specific setting of model (\ref{model}). Lower bounds
for minimax risk are classically derived by controlling the probability
for the likelihood ratio (in the statistical model of interested) of
being strictly greater than some constant uniformly over an appropriate
set of test functions. To derive Theorem \ref{theo:lowerbound}, we
show that one needs to control the expectation over the random shifts
of the likelihood ratio associated to model (\ref{model}), and not
only the likelihood ratio itself. Hence, the proof of Theorem \ref
{theo:lowerbound} is not a direct and straightforward adaptation of
Assouad's cube technique or Lemma 10.1 in \cite{HKPT} as used
classically in a standard white noise model to derive minimax risk in
nonparametric deconvolution in the dense case. For more details, we
refer to the proof of Theorem \ref{theo:lowerbound} in the \hyperref[app]{Appendix}.

Deriving minimax risk in the dense case for the model (\ref{model}) is
rather difficult and the proof is quite long and technical. In the
sparse case, finding lower bounds for the minimax rate is also a
difficult task. We believe that this could be done by adapting to model
(\ref{model}) a result by \cite{korotsy} which yields a lower bound
for a specific problem of distinguishing between a finite number of
hypotheses (see Lemma 10.1 in \cite{HKPT}). However, this is far
beyond the scope of this paper and we leave this problem open for
future wok.

\section{Estimating $f$ when the density $g$ is unknown} \label
{sec:estimgamma}

Obviously, assuming that the density $g$ of the shifts is known is not
very realistic in practice. However, estimating $f$ when the density
$g$ is unknown falls into the setting of inverse problems with an
unknown operator which is a difficult problem. Recently, some papers
\cite{cavraim,EfromovichK01,hoffreiss,neu97}
have considered nonparametric estimator for inverse problem
with a partially unknown operator, by assuming that an independent
sample of noisy eigenvalues is available which allows to build an
estimator of the $\gamma_{\ell}$'s. In the settings of these papers,
the distribution of the noisy eigenvalues sample is supposed to be
known (typically Gaussian). However, in model (\ref{model}), such
assumptions are not realistic, and therefore a data-based estimator of
$g$ has to be found. For this purpose, we propose to make a connection
between mean pattern estimation in model (\ref{model}) and well-known
results on Frechet mean estimation for manifold-valued variables; see,
for example, \cite{batach1,batach2}.

\subsection{Frechet mean for functional data} \label{sec:frechet}

Suppose that $Z_{1},\ldots,Z_{n}$ denote i.i.d. random variables
taking their values in a vector space $V$. As $V$ is a linear space
(with addition well defined), an estimator of a mean pattern for the
$Z_{m}$'s is given by the usual linear average $\overline{Z}_{n} =
\frac{1}{n} \sum_{m=1}^{n} Z_{m}$. However, in many applications, some
geometric and statistical considerations may lead to the assumption
that two vectors $Z,Z'$ in $V$ are considered to be the same if they
are equal up to certain transformations which are represented by the
action of some group $H$ on the space~$V$. A well-known example (see
\cite{batach1,batach2} and references therein) is the case
where $V = {\mathbb{R}}^{2 \times k}$, the space of $k$ points in the plane
${\mathbb{R}}^{2}$, and $H$ is generated by composition of scaling, rotations
and translations of the plane, namely
\[
h \cdot Z = a \pmatrix{
\cos\theta& - \sin\theta\cr\sin\theta& \cos
\theta
} Z + b,
\]
for $h = h(a,\theta,b) \in H$, with $(a,\theta,b) \in{\mathbb
{R}}^{+} \times
[0, 2\pi] \times{\mathbb{R}}^{2}$. In this setting, two vectors
$Z,Z' \in{\mathbb{R}}
^{2 \times k}$ represent the same shape if
\[
d_{H}(Z,Z') := {\inf_{(a,\theta,b) \in{\mathbb{R}}^{+} \times[0,
2\pi]
\times{\mathbb{R}}^{2}}} \|Z - h(a,\theta,b) \cdot Z' \|_{{\mathbb
{R}}^{2k}} = 0,
\]
which leads to Kendall's shape space $\Sigma_{2}^{k}$ consisting
of the equivalent classes of shapes in ${\mathbb{R}}^{2 \times k}$
under the
action of scaling, rotations and translations (see, e.g.,
\cite{batach1,batach2} and references therein). Since the space
$\Sigma_{2}^{k}$ is a nonlinear manifold, the usual linear average
$\overline{Z}_{n}$ does not fall into $\Sigma_{2}^{k}$ due to the
fact that the Euclidean distance \mbox{$\| \cdot\|_{{\mathbb{R}}^{2 \times
k}}$} is
not meaningful to represent shape variations. A better notion of
empirical mean $\tilde{Z}_{n}$ of $n$ shapes in ${\mathbb{R}}^{2
\times k}$ is
given by (see, e.g., \cite{batach1}):
$\tilde{Z}_{n} = \arg\min_{Z \in\Sigma_{2}^{k}} \frac{1}{n} \sum
_{m=1}^{n} d^{2}_{H}(Z,Z_{m})$.
More generally, Frechet \cite{frech} has extended the notion of
averaging to general metric spaces via mean squared error minimization
in the following way: if $Z_{1},\ldots,Z_{n}$ are i.i.d. random
variables in a general metric space $\mathcal M$, with a distance $d
\dvtx\mathcal M
\times\mathcal M\to{\mathbb{R}}^{+}$, then the Frechet mean of a
collection of such
data points is defined as the minimizer (not necessarily unique) of the
sum-of-squared distances to each of the data points, that is
\[
\tilde{Z}_{n} = \mathop{\arg\min}_{Z \in\mathcal M} \frac{1}{n} \sum_{m=1}^{n}
d^{2}(Z,Z_{m}).
\]
Now let us return to the randomly shifted curve model (\ref{model}).
Define $H = {\mathbb{R}}$ as the translation group acting on periodic functions
$f \in L^{2}([0,1])$ with period 1 by
\[
\tau\cdot f (x) = f (x + \tau)\qquad \mbox{for } x \in[0,1] \mbox
{ and } \tau\in H.
\]
Let $Y_{1},\ldots,Y_{n}$ be $n$ functions (possibly random) in
$L^{2}([0,1])$. Following the definition of Frechet mean, a notion of
average for functional data taking into account the action of the
translation group $H = {\mathbb{R}}$ would be
\[
\tilde{f}_{n} = \mathop{\arg\min}_{f \in L^{2}([0,1])} \frac{1}{n} \sum
_{m=1}^{n} \min_{\tau_{m} \in{\mathbb{R}}} \int_{0}^{1} {|f(x) -
Y_{m}(x+\tau_{m})|}^{2}\,dx.
\]
If the $Y_{m}$'s are noisy curves generated from the randomly shifted
curve model (\ref{model}), a presmoothing step of the observed curves
seems natural to compute a consistent Frechet mean estimate. In the
case of the translation group, this smoothing step and the definition
of Frechet mean can be expressed in the Fourier domain as
%
%
\begin{equation}\qquad
(\hat{\theta}_{-\ell_{0}},\ldots,\hat{\theta}_{\ell_{0}}) =
{\mathop{\arg\min}_{(\theta_{-\ell_{0}},\ldots,\theta_{\ell_{0}}) \in
{\mathbb{R}}
^{2 \ell_{0} + 1}} \frac{1}{n} \sum_{m=1}^{n} \min_{\tau_{m} \in
{\mathbb{R}}} \sum_{|\ell| \leq\ell_{0} }} {|c_{m,\ell} e^{2 i \ell
\pi\tau
_{m}} - \theta_{\ell} |^{2}},
\end{equation}
where $c_{m,\ell} = \int_{0}^1 e^{- 2 i \ell\pi x} \,dY_{m}(x)$ and
$\ell_{0}$ is some frequency cut-off parameter whose choice will be
discussed later. A smoothed Frechet mean is then given by
$\tilde{f}_{n,\ell_{0}} = \sum_{|\ell| \leq\ell_{0} } \hat{\theta
}_{\ell} e^{- 2 i \ell\pi x}$.
Then define the following criterion for $\bolds\tau= (\tau
_{1},\ldots
,\tau_{n}) \in{\mathbb{R}}^{n} $:
%
%
\begin{eqnarray}
M_{n}(\bolds\tau) & = & \frac{1}{n} \sum_{m=1}^{n} \sum_{|\ell|
\leq
\ell_{0} } \Biggl|c_{m,\ell} e^{2 i \ell\pi\tau_{m}} - \frac
{1}{n} \Biggl( \sum_{q=1}^{n} c_{q,\ell} e^{2 i \ell\pi\tau_{q}}
\Biggr) \Biggr|^{2},
\end{eqnarray}
and remark that the computation of $\overline{f}_{n,\ell_{0}}$ can be made
in two steps since it can be checked that
$\hat{\theta}_{\ell} = \frac{1}{n} \sum_{m=1}^{n} c_{m,\ell} e^{2
i \ell\pi\hat{\tau}_{m}}$,
where
%
%
\begin{equation} \label{eq:critFrechet}
(\hat{\tau}_{1},\ldots,\hat{\tau}_{n}) = \mathop{\arg\min}_{ (\tau
_{1},\ldots,\tau_{n}) \in{\mathbb{R}}^{n} } M_{n}(\tau_{1},\ldots
,\tau_{n}).
\end{equation}
Therefore, computing the Frechet mean of the smoothed curves
$Y_{1},\ldots,Y_{n}$ requires minimization of the above criteria which
automatically yields estimators of the random shifts $\tau_{1},\ldots
,\tau_{n}$ in model (\ref{model}). This allows us to construct an
estimator of the common shape $f$ by $\tilde{f}_{n,\ell_{0}}$ in the
case of an unknown density $g$, and the estimates $(\hat{\tau
}_{1},\ldots,\hat{\tau}_{n})$ of the random shifts can be used to
estimate the density $g$ itself and the eigenvalues $\gamma_{\ell}$.
The goal of this section is thus to study some statistical properties
of such a two-step procedure which, to the best of our knowledge, has
not been considered before in the setting of model (\ref{model}) and
in connection with Frechet mean for functional data. Moreover, it will
be shown in our numerical experiments that the criterion (\ref
{eq:critFrechet}) can be minimized using a gradient-descent algorithm
which leads to a new and fast method for estimating $f$ in the case of
an unknown density~$g$.

\subsection{Upper bound for the estimation of the shifts}

Recall that our model (\ref{model}) in the Fourier domain is
%
%
\begin{equation} \label{eq:mod2}
c_{m,\ell} = \theta_{\ell} e^{-i 2 \pi\ell\tau_{m}^{\ast}} +
\varepsilon z_{\ell,m} ,\qquad \ell\in{\mathbb{Z}}\mbox{ for }
m=1,\ldots,n,
\end{equation}
where $z_{\ell,m}$ are i.i.d. $N_{{\mathbb{C}}} (0,1 )$ variables,
the random shifts $\tau_{m}^{\ast},m=1,\ldots,n$, are i.i.d.
variables with density $g$, and $\theta_{\ell} = \int_{0}^{1} f(x)
e^{-i 2 \pi\ell x}\,dx $. Model (\ref{eq:mod2}) is clearly
nonidentifiable, as for any $\tau_{0} \in{\mathbb{R}}$, one can
replace the
$\theta_{\ell}$'s by $\theta_{\ell}e^{i 2 \pi\ell\tau_{0}}$ and
the $\tau_{m}^{\ast}$'s by $\tau_{m}^{\ast}-\tau_{0}$ without
changing the formulation of the model. Let us thus introduce the
following identifiability conditions.
\begin{ass} \label{ass:g}
The density $g$ has a compact support included in the interval
${\mathcal T}=
[-\frac{1}{4},\frac{1}{4}]$ and has zero mean, that is, is
such that
$\int_{{\mathcal T}} \tau g(\tau) \,d \tau= 0$.
\end{ass}
\begin{ass} \label{ass:indent}
The unknown shape function $f$ is such that $\theta_{1} \neq0$.
\end{ass}

Let
$\overline{{\mathcal T}}_{n} = \{ (\tau_{1},\ldots,\tau_{n}) \in
{\mathcal T}^{n}
\mbox{ such that } \sum_{m=1}^{n} \tau_{m} = 0 \}$.
Using the identifiability condition given by Assumption
\ref{ass:g},
it is natural to define estimators of the true shifts $\tau_{1}^{\ast
},\ldots,\tau_{n}^{\ast}$ as
\[
\hat{\bolds\tau} = (\hat{\tau}_{1},\ldots,\hat{\tau}_{n}) =
\mathop{\arg\min}
_{ \bolds\tau\in\overline{{\mathcal T}}_{n} } M_{n}(\bolds\tau),
\]
that is, by considering the estimators that minimize the empirical
criterion $M_{n}(\bolds\tau)$ on the constrained set $ \overline
{{\mathcal T}
}_{n}$. Then the following theorem holds.
\begin{theo} \label{th:estimshift}
Suppose that Assumptions \ref{ass:g} and \ref{ass:indent} hold. Then
for any $t > 0$
%
%
\begin{equation} \label{eq:cons}
{\mathbb P} \Biggl( \frac{1}{n} \sum_{m=2}^{n} (\hat{\tau}_{m} -\tau
_{m}^{\ast} )^{2} \geq C(f,\ell_{0},\varepsilon,n,t,g) \Biggr) \leq3
\exp(-t)
\end{equation}
with $C(f,\ell_{0},\varepsilon,n,t,g) = 4 \max[ C_{1}(f,\ell
_{0}) ( \sqrt{C_{2}(\varepsilon,n,\ell_{0},t)} + C_{2}(\varepsilon
,n,\ell_{0},t) )$,\break $C_{3}(t,n,g) ] $, where $C_{1}(f,\ell
_{0})$ is a positive constant depending only on the shape function $f$
and the frequency cut-off parameter $\ell_{0}$,
\[
C_{2}(\varepsilon,n,\ell_{0},t) = \varepsilon^2 (2 \ell_{0} + 1) + 2
\varepsilon^2 \sqrt{\frac{2 \ell_{0} + 1}{n}t} + 2 \frac{\varepsilon^2}{n}t
\]
and
\[
C_{3}(t,n,g) = \Biggl( \sqrt{2 \sigma^{2}_{g} \frac{t}{n} } + \frac
{t}{12 n} \Biggr)^{2} \qquad\mbox{with } \sigma^{2}_{g} = \int
_{{\mathcal T}}
\tau^{2} g(\tau) \,d \tau.
\]
\end{theo}

Theorem \ref{th:estimshift} provides an upper bound (in probability)
for the consistency of the estimators $\hat{\tau}_{m}$ of the true
random shifts $\tau_{m}^{\ast},m=2,\ldots,n$, using the standard
squared distance. Note that since the minimum of $M_{n}(\bolds\tau)$ is
computed on the constrained set $ \overline{{\mathcal T}}_{n}$, it follows
that $\hat{\tau}_{1} = - \sum_{m=2}^{n} \hat{\tau}_{m}$. However,
one can remark that as $n \to+ \infty$, the constant $C(f,\ell
_{0},\varepsilon,n,t,g)$ in inequality (\ref{eq:cons}) tends to $4
C_{1}(f,\ell_{0}) ( \varepsilon^2 (2 \ell_{0} + 1) + \varepsilon
\sqrt{ 2 \ell_{0} + 1} )$. This shows that $\hat{\tau}_{m},
m=2,\ldots,n$, are not consistent estimators in the sense that
inequality (\ref{eq:cons}) cannot be used to prove that $\lim_{n \to
+ \infty} \frac{1}{n} \sum_{m=2}^{n} (\hat{\tau}_{m} -\tau
_{m}^{\ast} )^{2} = 0$ in probability. On the contrary, inequality
(\ref{eq:cons}) suggests that there exists a constant $C >0$ such that
$\frac{1}{n} \sum_{m=2}^{n} (\hat{\tau}_{m} -\tau_{m}^{\ast}
)^{2} > C ( \varepsilon^2 (2 \ell_{0} + 1) + \varepsilon\sqrt{ 2
\ell_{0} + 1} )$ with positive probability, and that the
accuracy of the estimates $\hat{\tau}_{m}, m=2,\ldots,n$, should
depend on the level of noise $\varepsilon^{2}$ and the frequency cut-off
$\ell_{0}$.

The choice of the frequency cut-off $\ell_{0}$ used to compute these
estimators is a delicate model selection problem. Theorem \ref
{th:estimshift} suggests that this choice should depend on $n$ and the
level of noise $\varepsilon$, but finding data-based values for $\ell
_{0}$ remains a challenge that we leave open for future work.

\subsection{Lower bound for the estimation of the shifts}

Let us now prove that the consistency of any estimate of the random
shifts in model (\ref{eq:mod2}) is limited by the level of noise
$\varepsilon^{2}$ in the observed curves. For this, let us make the
following smoothness assumptions.
\begin{ass} \label{ass:f}
The function $f$ is such that $\sum_{\ell\in{\mathbb{Z}}} (2\pi
\ell)^{2}
|\theta_{\ell}|^{2} < + \infty$.
\end{ass}
\begin{ass} \label{ass:gT}
The density $g$ is compactly supported on a interval ${\mathcal T}=
[\tau
_{\min},\tau_{\max}]$ such that $\lim_{\tau\to\tau_{\min}}
g(\tau) = \lim_{\tau\to\tau_{\max}} g(\tau) = 0$.
\end{ass}

Then, using general results on the Van Tree's inequality \cite{gillev}
in model (\ref{eq:mod2}), the following theorem holds.
\begin{theo} \label{theo:VanTree}
Denote by $X = (c_{m,\ell})_{\ell\in{\mathbb{Z}}, m=1,\ldots, n}$
the set
of observations taking values in the set ${\mathcal X}= {\mathbb
{C}}^{\infty\times
n}$. Let $\hat{\tau}^{n} = \hat{\tau}^{n}(X)$ denote any estimator
(a measurable function of the observations $X$) of the true shifts
$(\tau_{1},\ldots,\tau_{n})$. Then, under Assumptions \ref{ass:f}
and \ref{ass:gT},
\[
{\mathbb E} \Biggl( \frac{1}{n} \sum_{m=1}^{n} (\hat{\tau}^{n}_{m}-
\tau
_{m}^{\ast})^{2} \Biggr) \geq\frac{\varepsilon^{2} }{\sum_{\ell\in
{\mathbb{Z}}} (2\pi\ell)^{2} |\theta_{\ell}|^{2} + \varepsilon^{2}
\int_{{\mathcal T}
} ( {\partial}/{\partial\tau} \log g(\tau) )^{2}
g(\tau) \,d \tau}.
\]
\end{theo}

Clearly, Theorem \ref{theo:VanTree} shows that as $n \to+ \infty$
then $ {\mathbb E} ( \frac{1}{n} \sum_{m=1}^{n} (\hat{\tau}^{n}_{m}-
\tau_{m}^{\ast})^{2} ) $ does not converge to zero which
explains the results obtained in Theorem \ref{th:estimshift} on the
consistency of the estimators $\hat{\tau}_{m}, m=2,\ldots,n$, based
on Frechet mean for functional data. Note that Assumption \ref{ass:f}
can be avoided if one only considers estimators $\hat{\tau}^{n,\ell
_{0}}$ of the shifts based on the observations $c_{m,\ell}$ for
$m=1,\ldots,n$ and $|\ell| \leq\ell_{0}$ in model~(\ref{eq:mod2}).
In this case, the lower bound in Theorem \ref{theo:VanTree} becomes
\[
{\mathbb E} \Biggl( \frac{1}{n} \sum_{m=1}^{n} (\hat{\tau}^{n,\ell
_{0}}_{m}- \tau_{m}^{\ast})^{2} \Biggr) \geq\frac{\varepsilon^{2}
}{\sum_{|\ell| \leq\ell_{0}} (2\pi\ell)^{2} |\theta_{\ell}|^{2}
+ \varepsilon^{2} \int_{{\mathcal T}} ( {\partial}/{\partial
\tau}
\log g(\tau) )^{2} g(\tau) \,d \tau}.
\]

\subsection{Estimation of the mean pattern $f$ and the density $g$}

An estimator of the eigenvalue $\gamma_{\ell}$ is given by
%
%
\begin{equation}\label{ea:hatgamma}
\hat{\gamma}_{\ell} = \frac{1}{n} \sum_{m=2}^{n} e^{- i 2 \pi\ell
\hat{\tau}_{m}}
\end{equation}
for $|\ell| \leq\ell_{0}$ and an estimator for the density $g$ is
naturally given by
$\hat{g}(x) = \sum_{|\ell| \leq\ell_{0}} \hat{\gamma}_{\ell}
e^{- i 2 \pi\ell x}$.
The mean pattern $f$ can be estimated by the smoothed Frechet mean
$\tilde{f}_{n,\ell_{0}}$ defined in Section \ref{sec:frechet}, but
following the results in Section \ref{sec:estimfnonlinear} on
nonlinear wavelet-based estimation, two other estimators for $f$ can be
defined: the first one is given by
%
%
\begin{equation} \label{eq:hatfnonlin1}
\hat{f}_{n,1} = \sum_{k =0}^{2^{j_{0}}-1} \hat{c}_{j_{0},k,1} \phi
_{j_{0},k} + \sum_{j = j_{0}}^{j_{1}} \sum_{k=0}^{2^{j}-1} \hat
{\beta}_{j,k,1} \mathbh{1}_{\{| \hat{\beta}_{j,k,1} | \geq
\hat
{\lambda}_{j,1}\}} \psi_{j,k},
\end{equation}
where $ \hat{\beta}_{j,k,1} = \sum_{\ell\in\Omega_{j}} \psi
^{j,k}_{\ell} \hat{\theta}_{\ell,1}$ and $\hat
{c}_{j_{0},k,1} = \sum_{\ell\in\Omega_{j_{0}}} \phi
^{j_{0},k}_{\ell} \hat{\theta}_{\ell,1}$ with
$\hat{\theta}_{\ell,1} = \frac{1}{\hat{\gamma}_{\ell}}\times\break (
\frac{1}{n} \sum_{m=1}^n c_{\ell,m} )$,
and
\[
\hat{\lambda}_{j,1} = 2 \hat{\sigma}_{j} \sqrt{\frac{2
\eta\log(n)}{n}}
\]
is the threshold suggested by the expression (\ref{eq:hatthr}) of
$\lambda_{j}$ with
$\hat{\sigma}^{2}_{j} = 2^{-j} \varepsilon^{2} \times\break\sum_{\ell\in\Omega
_{j}} | \hat{\gamma}_{\ell} |^{-2}$.
A second estimator is given by first realigning the curves using the
estimators of the shifts, namely,
%
%
\begin{equation} \label{eq:hatfnonlin2}
\hat{f}_{n,2} = \sum_{k =0}^{2^{j_{0}}-1} \hat{c}_{j_{0},k,2} \phi
_{j_{0},k} + \sum_{j = j_{0}}^{j_{1}} \sum_{k=0}^{2^{j}-1} \hat
{\beta}_{j,k,2} \mathbh{1}_{\{| \hat{\beta}_{j,k,2} | \geq
\hat
{\lambda}_{j,2}\}} \psi_{j,k},
\end{equation}
where $ \hat{\beta}_{j,k,2} = \sum_{\ell\in\Omega_{j}} \psi
^{j,k}_{\ell} \hat{\theta}_{\ell,2}$ and $\hat
{c}_{j_{0},k,2} = \sum_{\ell\in\Omega_{j_{0}}} \phi
^{j_{0},k}_{\ell} \hat{\theta}_{\ell,2}$ with
\[
\hat{\theta}_{\ell,2} = \frac{1}{n} \sum_{m=2}^n c_{\ell,m} e^{i
2 \pi\ell\hat{\tau}_{m}},
\]
and $\hat{\lambda}_{j,2}$ is a threshold whose choice would depend on
the law of the $\hat{\beta}_{j,k,2}$'s. Studying the consistency and
the rate of convergence of the estimators $\hat{f}_{n,1}$ and $\hat
{f}_{n,2}$ is a difficult task. Indeed the results in Section \ref
{sec:estimfnonlinear} have been derived using the fact that the law of
the wavelet coefficients $\hat{\beta}_{j,k}$ and $\hat{c}_{j_{0},k}$
given by (\ref{eq:empwavcoef}) is \textit{known} which allows the
calibration of the threshold $\lambda_{j}$ in (\ref{eq:hatthr}).
Thus, we simply suggest to take $\hat{\lambda}_{j,2} = \hat{\lambda}_{j,1}$.
Extending the asymptotic results of Section \ref{sec:estimfnonlinear}
remains a challenge that is beyond the scope of this paper. Moreover,
the results of Theorems~\ref{th:estimshift} and \ref{theo:VanTree}
suggest that the estimators $\hat{f}_{n,1}$ and $\hat{f}_{n,2}$ could
be consistent by considering a double asymptotic setting with $n \to+
\infty$ and $\varepsilon\to0$ which is an interesting point of view for
future work that certainly leads to different minimax rates of convergence.

\section{Numerical experiments} \label{sec:simu}

We compare our approach with the Procrustean mean which is a standard
algorithm commonly used to extract a mean pattern. The Procrustean mean
is based on an alternative scheme between estimation of the shifts and
averaging of back-transformed curves given estimated values of the
shifts parameters; see, for example, \cite{wanggas,kg}. To be
more precise, it consists of an initialization step $\hat{f}_{0} =
\frac{1}{n} \sum_{m=1}^{n} Y_{m}$ which is the simple average of the
observed curves, that is taken as a first reference mean pattern. Then,
at iteration $1 \leq i \leq i_{\max}$, it computes for all $1 \leq m
\leq n$ the estimators $\hat{\tau}_{m,i}$ of the $m$th shift as $\hat
{\tau}_{m,i} = \arg\min_{\tau\in{\mathbb{R}}} \|Y_{m}(\cdot+\tau
) - \hat
{f}_{i-1} \|^{2} $ and then takes $ \hat{f}_{i}(x) = \frac{1}{n} \sum
_{m=1}^{n} Y_{m}(x +\hat{\tau}_{m,i}) $ as a new reference mean
pattern. This is repeated until the estimated reference curve does not
change, and usually the algorithm converges in a few steps (we took
$i_{\max}=3$). In all simulations, we used the wavelet toolbox Wavelab
\cite{BCDJ95} and the WaveD algorithm developed by \cite{rast07} for
fast deconvolution with Meyer wavelets.

\subsection{Shift estimation by gradient descent} \label{sec:grad}

Let us denote by $\nabla M_{n} (\bolds\tau) \in{\mathbb{R}}^{n}$
the gradient of
$M_{n}(\bolds\tau)$ at $\bolds\tau\in{\mathbb{R}}^{n}$. This
gradient is simple to
compute as for $m=1,\ldots,n$:
\[
\frac{\partial}{\partial\tau_{m}} M_{n}(\bolds\tau) = -\frac{2}{n}
\sum_{|\ell| \leq\ell_{0} } \Re\Biggl[ 2 i \pi\ell c_{\ell
,m}e^{2 i \ell\pi\tau_{m}} \Biggl( \overline{ \frac{1}{n} \sum
_{q=1}^{n} c_{\ell,q} e^{2 i \ell\pi\tau_{q}} } \Biggr) \Biggr].
\]
In practice, to estimate the shifts, the criterion $M_{n}(\bolds\tau
)$ is
minimized by the following gradient descent algorithm with the
constraint that $\tau_{1} = - \sum_{m=2}^{n} \tau_{m}$:

\textit{Initialization}: let $\bolds\tau^{0} = 0 \in
{\mathbb{R}}^{n}$,
$\delta_{0} = \frac{1}{\| \nabla M_{n} (\bolds\tau^{0}) \|}$, $M(0)
=M_{n}(\bolds\tau^{0})$ and set $p = 0$.\vspace*{1pt}

\textit{Step} 2: let $\bolds\tau^{\mathrm{new}} = \bolds\tau^{p} -
\delta_{p}
\nabla M_{n}(\bolds\tau^{p})$ and
$\tau^{\mathrm{new}}_{1} = - \sum_{m=2}^{n} \tau^{\mathrm{new}}_{m}$.

Let $M(p+1) =M_{n}(\bolds\tau^{\mathrm{new}})$.

\textit{While $M(p+1) > M(p)$ do}
\[
\delta_{p} = \delta_{p} / \kappa\quad \mbox{and}\quad \bolds\tau
^{\mathrm{new}} =
\bolds\tau^{p} - \delta_{p} \nabla M_{n}(\bolds\tau^{m})\qquad \mbox{with } \tau
^{\mathrm{new}}_{1} = - \sum_{m=2}^{n} \tau^{\mathrm{new}}_{m},
\]
and set $M(p+1) =M_{n}(\bolds\tau^{\mathrm{new}})$.

\textit{End while}

Then, take $\bolds\tau^{p+1} = \bolds\tau^{\mathrm{new}}$.

\textit{Step} 3: if $M(p) - M(p+1) \geq\rho(M(1) - M(p+1))$
then set $p = p + 1$ and return to Step 2, else stop the iterations,
and take $\hat{\bolds\tau} = \bolds\tau^{p+1}$.

In the above algorithm, $\rho> 0$ is a small stopping parameter and
$\kappa> 1$ is a parameter to control the choice of the adaptive step
$\delta_{p}$.

\subsection{Estimation with an unknown density $g$}

\begin{figure}

\includegraphics{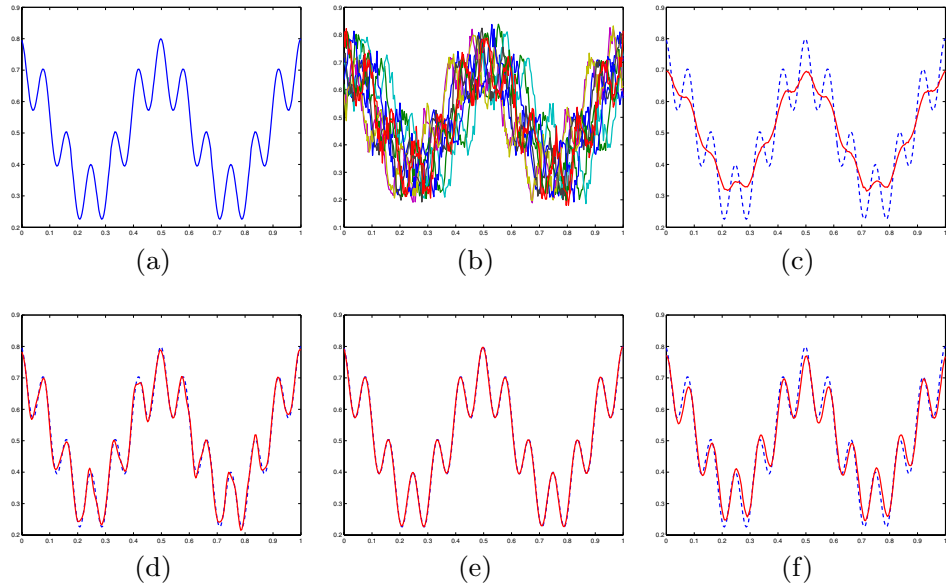}

\caption{Wave function. \textup{(a)} Mean pattern $f$, \textup{(b)} sample of 10 curves
out of $n=200$, \textup{(c)} direct mean, deconvolution by wavelet thresholding
with \textup{(d)} $\hat{f}_{n,1}$ and \textup{(e)} $\hat{f}_{n,2}$, \textup{(f)}
Procrustean mean.}
\label{Fig:test1unknown}
\end{figure}

\begin{figure}

\includegraphics{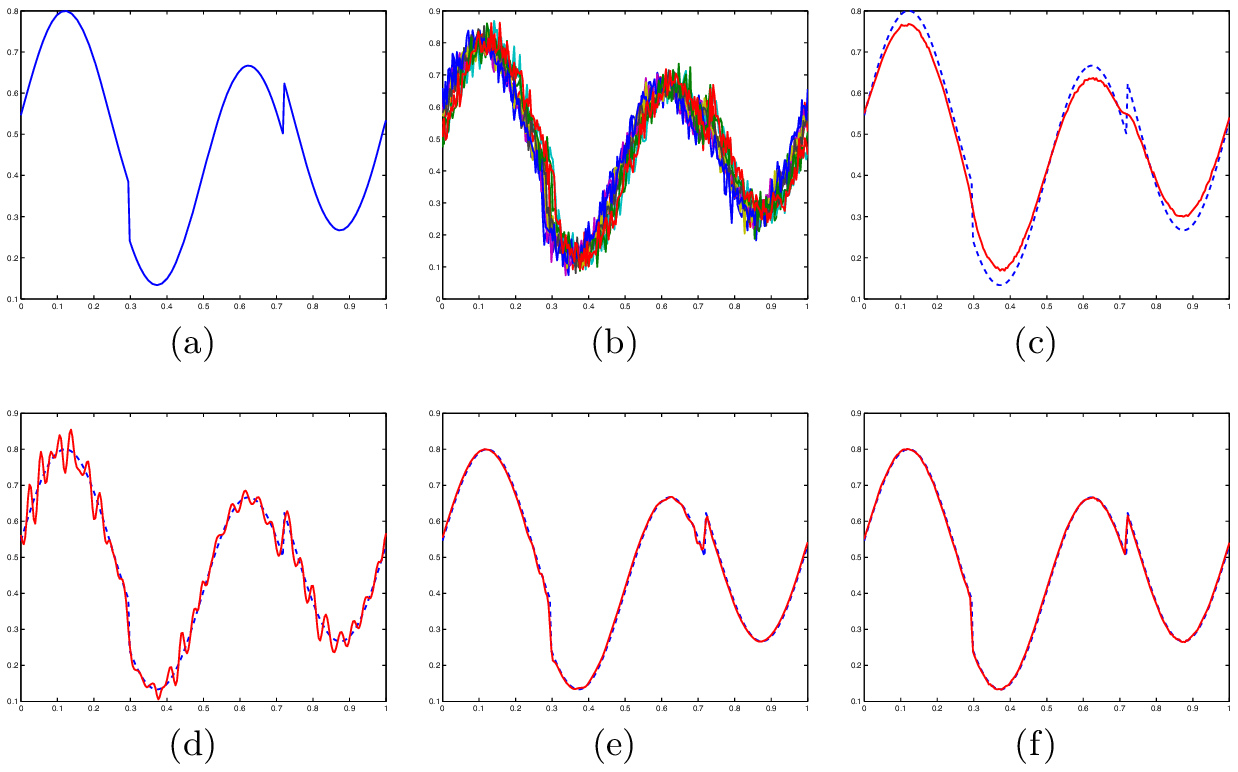}

\caption{HeaviSine function. \textup{(a)} Mean pattern $f$, \textup{(b)} sample of 10
curves out of $n=200$, \textup{(c)} direct mean, deconvolution by wavelet
thresholding with \textup{(d)} $\hat{f}_{n,1}$ and \textup{(e)} $\hat{f}_{n,2}$,
\textup{(f)} Procrustean mean.} \label{Fig:test2unknown}
\end{figure}

\begin{figure}

\includegraphics{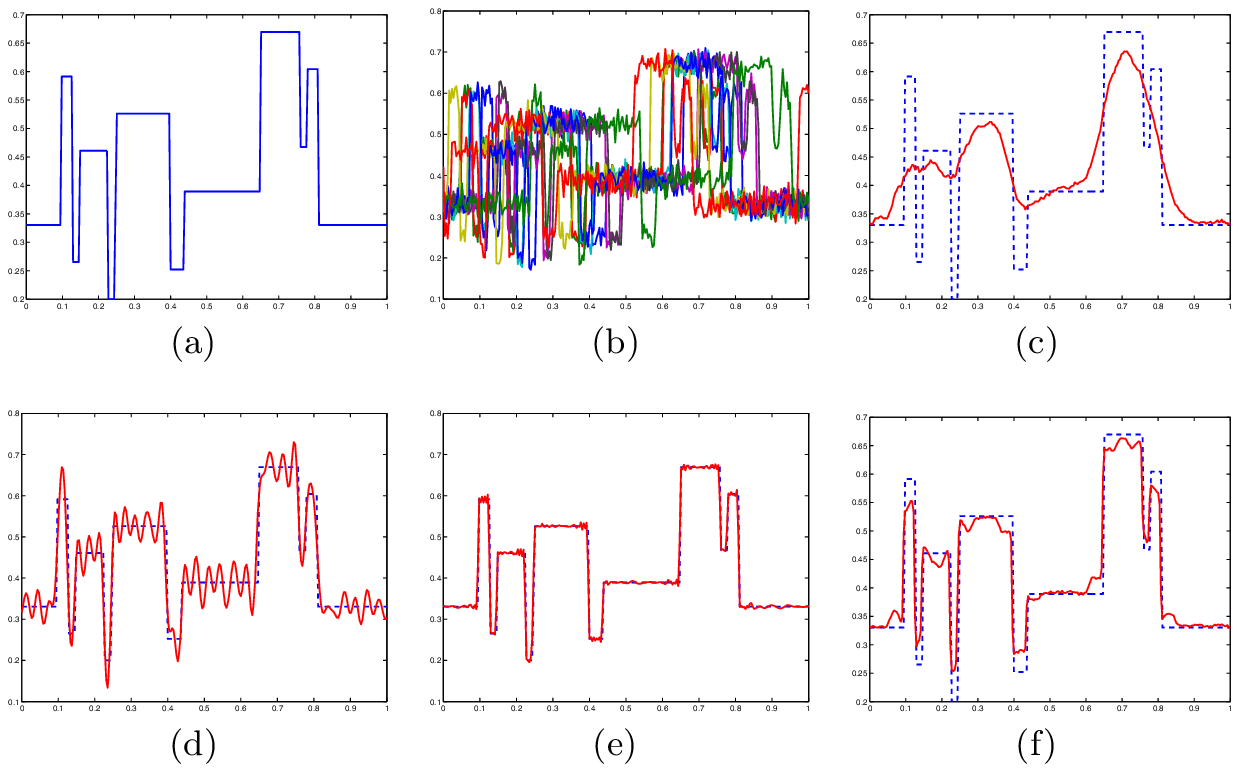}

\caption{Blocks function. \textup{(a)} Mean pattern $f$, \textup{(b)} sample of 10
curves out of $n=200$, \textup{(c)} direct mean, deconvolution by wavelet
thresholding with \textup{(d)} $\hat{f}_{n,1}$ and \textup{(e)} $\hat{f}_{n,2}$,
\textup{(f)} Procrustean mean.} \label{Fig:test3unknown}
\end{figure}

\begin{figure}

\includegraphics{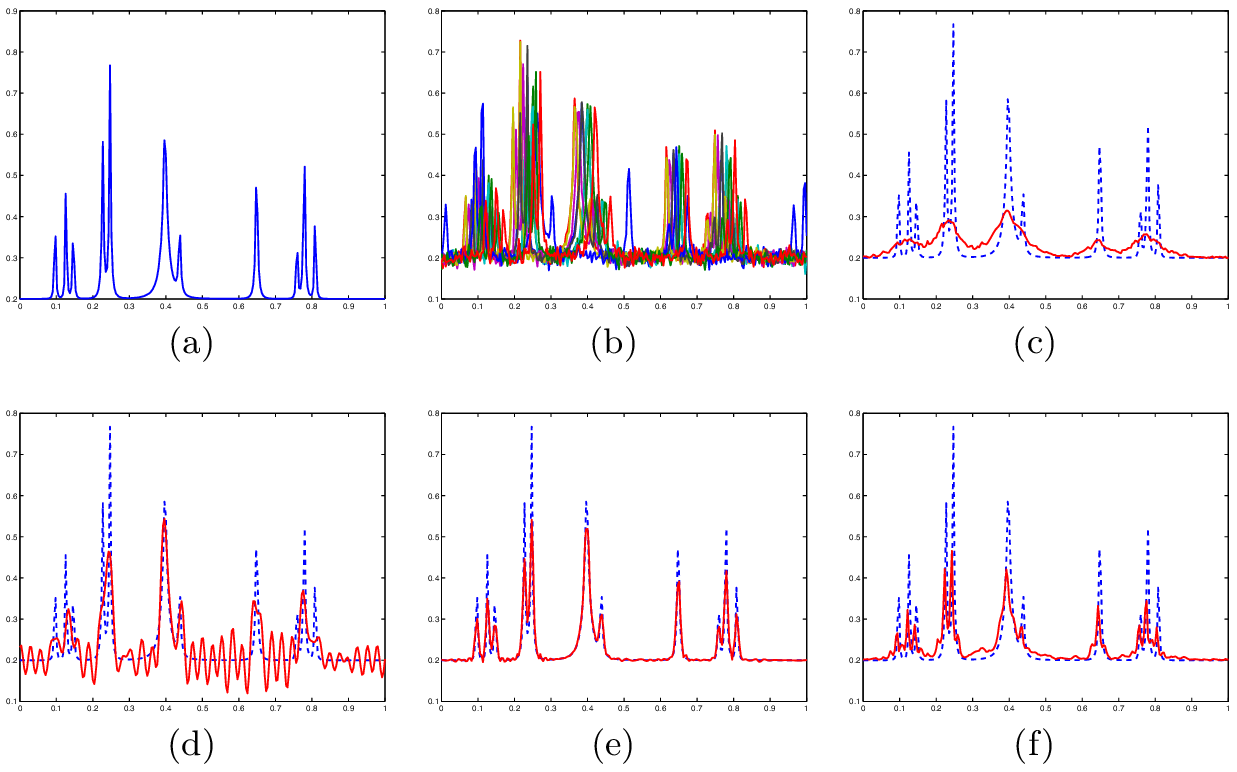}

\caption{Bumps function. \textup{(a)} Mean pattern $f$, \textup{(b)} sample of 10 curves
out of $n=200$, \textup{(c)} direct mean, deconvolution by wavelet thresholding
with \textup{(d)} $\hat{f}_{n,1}$ and \textup{(e)} $\hat{f}_{n,2}$, \textup{(f)} Procrustean mean.}
\label{Fig:test4unknown}
\end{figure}

For the mean pattern $f$ to recover, we consider the four tests
functions shown in Figures \ref{Fig:test1unknown}(a)--\ref
{Fig:test4unknown}(a). Then, we simulate $n = 200$ randomly shifted
curves with shifts following a Laplace distribution $g(x) = \frac
{1}{\sqrt{2}\sigma} \exp( -\sqrt{2}\frac{|x|}{\sigma}
)$ with $\sigma= 0.1$. Gaussian noise with a moderate variance
(different to that used in the Laplace distribution) is then added to
each curve. A~subsample of 10 curves is shown in Figures \ref
{Fig:test1unknown}(b)--\ref{Fig:test4unknown}(b) for each test function,
and the average of the observed curves, referred to as the direct mean
in what follows, is displayed in Figures \ref{Fig:test1unknown}(c)--\ref
{Fig:test4unknown}(c). Note this gives a poor estimator of the mean pattern.

%
%

%
%

%
%

%
%

The Fourier coefficients of the density $g$ are given by $\gamma_{\ell
} = \frac{1}{1+2\sigma^{2} \pi^{2} \ell^{2}}$ which corresponds to
a degree of ill-posedness $\nu=2$. An estimator $\hat{\gamma}_{\ell
}$ of $\gamma_{\ell}$ can be computed as explained in Section \ref
{sec:estimgamma} using the gradient descent algorithm described in
Section \ref{sec:grad}. In our simulations, we took the arbitrary
choice $\ell_{0} = 3$ for the frequency cut-off which gives
satisfactory results in the numerical experiments.

To compute the threshold, $\hat{\lambda}_{j,1} = \hat{\lambda
}_{j,2}$ used in the definition of $\hat{f}_{n,1}$ and $\hat
{f}_{n,2}$ (see Section \ref{sec:estimgamma}) one has to estimate
$\varepsilon^{2}$. This is done by taking $ \hat{\varepsilon}^{2} =
\frac
{1}{n} \sum_{m=1}^{n} \hat{\varepsilon}_{m}^{2} $, where the variance
$\varepsilon^{2}_{m}$ of the noise for the $m$th curve is easily
estimated using the wavelet coefficients at the finest resolution
level. Note that such thresholds are quite simple to compute using the
fast Fourier transform and the fact that the set of frequencies $
\Omega_{j}$ can be easily obtained using Wavelab. Finally, we have
found that choosing $\eta$ between 1 and 2 to compute $\hat{\lambda
}_{j,1}$ gives quite satisfactory results.

Then we took $j_{0} = 3 \approx\log_{2}( \log(n) )$, but the choice
$j_{1} \approx\frac{1}{2 \nu+1} \log_{2}( \frac{n}{\log(n)} )$ is
obviously too small. So in our simulations, $j_{1}$ is chosen to be the
maximum resolution level allowed by the discretization,
that is, $j_{1} = \log_2(N)-1 = 7$. For each test
function, the estimators $\hat{f}_{n,1}$, $\hat{f}_{n,2}$ are
displayed\vspace*{1pt} in Figures \ref{Fig:test1unknown}(d)(e)--\ref
{Fig:test4unknown}(d)(e). The Procrustean mean\vspace*{1pt} is displayed is Figures
\ref{Fig:test1unknown}(f)--\ref{Fig:test4unknown}(f). One can see
that the results are rather satisfactory for $\hat{f}_{n,1}$ and the
Procrustean mean. Clearly, the best results are given by the estimator
$\hat{f}_{n,2}$ which gives very good estimates of the function $f$
particularly for functions with isolated singularities such as the
Blocks and Bumps functions in Figures \ref{Fig:test3unknown} and \ref
{Fig:test4unknown}. It should be noted that these results are obtained
in the case of an unknown density $g$ which shows the quality of the
procedure proposed in Section \ref{sec:estimgamma} to estimate the
shifts and the $\gamma_{\ell}$'s. For reasons of space a detailed
simulation study is not given, but it has been found that the good
performances of the wavelet-based estimator remain consistent across
other standard test signals.

\section{Conclusions and future work}

This paper makes a connection between mean pattern estimation and the
statistical analysis of inverse problems for a very simple model with
shifted curves. A natural extension would be to consider more complex
deformations in SIM models such as the homothetic shifted regression
model proposed in \cite{vimond}, or the rigid deformation model for
images considered in~\cite{BGV}. The results obtained on the
nonconsistency of the estimation of the shifts can be elaborated in
future work on deformable models for signal and image analysis. In
particular, it would be interesting to obtain similar results for more
general deformation parameters in SIM models. Another interesting
question in such models is whether one can estimate an unknown mean
pattern consistently even if it is impossible to construct consistent
estimators for the deformation parameters such as the random shifts in
model (\ref{model}).

Another promising approach would be to consider a double asymptotic
setting with $n \to+ \infty$ and $\varepsilon\to0$ to study the
consistency and rate of convergence for estimators of the mean pattern
and the unknown density $g$.


\begin{appendix}\label{app}

\section*{Appendix}

In what follows, $C$ will denote a generic constant whose value may
change from line to line.
\begin{pf*}{Proof of Theorem \ref{theo:main}} It follows
immediately from Theorems \ref{theo:upperbound} and \ref
{theo:lowerbound}.
\end{pf*}
\begin{pf*}{Proof of Proposition \ref{prop:decomp}} Let $ \kappa
_\ell= (
\frac{\tilde{\gamma}_{\ell}}{\gamma_{\ell}}-1 ) \theta
_{\ell}$ and $ \varepsilon_{\ell,n} = \frac{\varepsilon}{\gamma_{\ell
}} ( \frac{1}{n} \sum_{m=1}^{n} z_{\ell,m} )$ for all
$\ell\in{\mathbb{Z}}$. Then, for a given filter $\lambda$, the risk
${\mathcal R}
(\hat{f}_{n,\lambda},f)$ can be written as
\begin{eqnarray*}
{\mathcal R}(\hat{f}_{n,\lambda},f) & = & \sum_{\ell
\in{\mathbb{Z}}}
(\lambda_{\ell} - 1)^2 \theta_{\ell} ^2 +
\mathbb{E} [\lambda_{\ell}^2 |\kappa_{\ell}|^2 + \lambda
_{\ell}^2 |\varepsilon_{\ell,n}|^2
]\\
&&{} +
\lambda_{\ell} (\lambda_{\ell}-1) \mathbb{E} [\overline{\theta}_{\ell
} \kappa_{\ell}+ \theta_{\ell} \overline{\kappa}_{\ell}]  +
\lambda_{\ell}(\lambda_{\ell}-1) \mathbb{E}[\theta_{\ell}\overline
{\varepsilon}_{\ell,n}+\overline{\theta}_{\ell}\varepsilon_{\ell,n}]\\
&&{} +
\lambda_{\ell}^2 \mathbb{E}[\overline{\kappa}_{\ell} \varepsilon_{\ell
,n}+\kappa_{\ell} \overline{\varepsilon}_{\ell,n}].
\end{eqnarray*}
Now using the fact that $\kappa_{\ell}$ and $\varepsilon_{\ell,n}$ are
independent and that ${\mathbb E}\varepsilon_{\ell,n} = 0$, we obtain that
\begin{eqnarray*}
R(\hat{f}_{n,\lambda},f) & = & \sum_{k\in\mathbb
{Z}} \biggl[ (\lambda_{\ell}-1)^2 |\theta_{\ell}|^2 + \lambda_{\ell
}^2 |\theta_{\ell}|^2
\mathbb{E} \biggl| \frac{\tilde{\gamma}_{\ell}}{\gamma_{\ell
}}-1 \biggr|^2 +
\frac{\lambda_{\ell}^2 \varepsilon^2}{n|\gamma_{\ell}|^2} \biggr]
\\
& =& \sum_{\ell\in{\mathbb{Z}}} (\lambda_{\ell
}-1)^2 |\theta
_{\ell}|^2 + \sum_{\ell\in{\mathbb{Z}}} \lambda
_{\ell}^2
|\theta_{\ell}|^2 \biggl( \frac{1}{|\gamma_{\ell}|^2} \biggl(
\frac{1}{n}+\frac{n-1}{n}
\gamma_{\ell} \gamma_{-\ell} \biggr) -1 \biggr) \\
& &{} + \sum_{\ell\in{\mathbb{Z}}}
\frac{\lambda_{\ell}^2 \varepsilon^2}{n |\gamma_{\ell}|^2} \\
& = &
\sum_{\ell\in{\mathbb{Z}}} (\lambda_{\ell}-1)^2
|\theta_{\ell}|^2 + \sum_{\ell\in{\mathbb{Z}}}
\frac{\lambda_{\ell}^2}{n} \biggl[|\theta_{\ell}|^2 \biggl(\frac
{1}{|\gamma_{\ell}|^2} -1 \biggr)
+ \frac{\varepsilon^2}{|\gamma_{\ell}|^2} \biggr],
\end{eqnarray*}
which completes the proof.
\end{pf*}
\begin{pf*}{Proof of Proposition \ref{prop:risklinearpoly}} From
Proposition \ref{prop:decomp} it follows that
\[
{\mathcal R}(\hat{f}_{n,\lambda^{M}},f) = \sum_{|\ell| > M_{n} }
|\theta
_{\ell}|^2 + \frac{1}{n} \sum_{|\ell| \leq M_{n}} \biggl(|\theta
_{\ell}|^2 \biggl(\frac{1}{|\gamma_{\ell}|^2} -1 \biggr) + \frac
{\varepsilon^2}{|\gamma_{\ell}|^2} \biggr).
\]
By assumption $f \in H_{s}(A)$, which implies that there exist two
positive constants $C_{1}$ and $C_{2}$ not depending on $f$ and $n$
such that for all sufficiently large $n$,
$\sum_{|\ell| > M_{n} } |\theta_{\ell}|^2 \leq C_{1} M_{n}^{-2s}$ and
$\frac{1}{n} \sum_{|\ell| \leq M_{n}} |\theta_{\ell}|^2 \leq C_{2} n^{-1}$.
Now, given that $g$ satisfies Assumption \ref{assordi}, it follows
that there exists a positive constant $C_{3}$ not depending on $f$ and
$n$ such that for all sufficiently large $n$,
$\frac{1}{n} \sum_{|\ell| \leq M_{n}} \frac{|\theta_{\ell}|^2 +
\varepsilon^2}{|\gamma_{\ell}|^2} \leq C_{3} n^{-1} M_{n}^{2 \nu+ 1}$.
Hence, the result immediately follows from the choice $M_{n} \sim
n^{{1}/({2s + 2 \nu+1})}$, which completes the proof.
\end{pf*}

For the proofs of Propositions \ref{prop:moments} and \ref
{prop:proba}, let us remark that that $\hat{\beta}_{j,k} - \beta
_{j,k} = Z_{1} + Z_{2}$ with
%
%
\begin{equation} \label{eq:Z1Z2}
Z_{1} = \sum_{\ell\in\Omega_{j}} \psi^{j,k}_{\ell} \theta_{\ell
} \biggl(\frac{\tilde{\gamma}_{\ell}}{ \gamma_{\ell}}-1\biggr) \quad\mbox{and}\quad
Z_{2} = \varepsilon\sum_{\ell\in\Omega_{j}} \frac{\psi
^{j,k}_{\ell}}{ \gamma_{\ell}} \eta_{\ell}.
\end{equation}

Under Assumption \ref{assdecay},
$G(x) = \sum_{m \in{\mathbb{Z}}} g(x+m)$
exists for all $x \in[0,1]$ and is a bounded density. Throughout the
proof, we use the following lemma whose proof is straightforward.
\begin{lemma} \label{lem:pars}
Let $h \in L^{2}([0,1])$ be a 1-periodic function on ${\mathbb{R}}$. Then,
$\int_{{\mathbb{R}}} h(x) \times\break g(x) \,dx = \int_{0}^{1} h(x) G(x) \,dx$.
\end{lemma}
%
%
\begin{pf*}{Proof of Proposition \ref{prop:moments}} First note
that since $|\psi^{j,k}_{\ell}| \leq2^{-j/2}$ and $\Omega_{j}
\subset[-2^{j+2}c_{0},-2^{j}c_{0}] \cup[2^{j}c_{0},2^{j+2}c_{0}]$
(see \cite{JKPR04jrssb}) it follows that $\#\{\Omega_{j}\} \leq4 \pi
2^{j}$ and that under Assumption \ref{assordi}, $|\gamma_{\ell
}|^{-2} \sim2^{2 j \nu}$ for all $\ell\in\Omega_{j}$. This implies
that there exists a constant $C > 0$ such that
%
%
\begin{equation} \label{eq:bound1}
\sum_{\ell\in\Omega_{j}} \biggl| \frac{\psi^{j,k}_{\ell}}{ \gamma
_{\ell}} \biggr|^{2} \leq C 2^{2 j \nu} \quad\mbox{and}\quad \sum
_{\ell\in\Omega_{j}} \biggl| \frac{\psi^{j,k}_{\ell}}{ \gamma_{\ell}}
\biggr| \leq C 2^{j (\nu+1/2)}.
\end{equation}
Then, we need the following lemma which shows that the Fourier
coefficients $\theta_{\ell} = \int_{0}^1 e^{- 2 i \ell\pi x} f(x)
\,dx$ are uniformly bounded for all $f \in B^{s}_{p,q}(A)$.
\begin{lemma} \label{lem:bound}
Let $1 \leq p \leq\infty$, $1 \leq q \leq\infty$, $s - 1/p +1/2 >
0$ and $A > 0$. Then there exists a constant $A' > 0$ such that for all
$f \in B^{s}_{p,q}(A)$ and all $\ell\in{\mathbb{Z}}$,
\mbox{$|\theta_{\ell}| \leq A'$.}
\end{lemma}
\begin{pf}
Since $|\phi^{j_{0},k}_{\ell}| \leq
2^{-j_{0}/2}$ and $|\psi^{j,k}_{\ell}| \leq2^{-j/2}$ one can remark
using the Cauchy--Schwarz inequality that for any $j_{0} \geq0$,
\begin{eqnarray*}
|\theta_{\ell}| & \leq& \sum_{k =0}^{2^{j_{0}}-1} |c_{j_{0},k}|
|\phi^{j_{0},k}_{\ell}| + \sum_{j = j_{0}}^{+ \infty} \sum
_{k=0}^{2^{j}-1} |\beta_{j,k}| |\psi^{j,k}_{\ell}| \\
& \leq& \Biggl( \sum_{k =0}^{2^{j_{0}}-1} |c_{j_{0},k}|^{2}
\Biggr)^{1/2} + \sum_{j = j_{0}}^{+ \infty} \Biggl( \sum_{k=0}^{2^{j}-1}
|\beta_{j,k}|^{2} \Biggr)^{1/2}.
\end{eqnarray*}
Now, using the inequality $ ( \sum_{r=1}^{m} |a_{r}|^{2}
)^{1/2} \leq m^{(1/2-1/p)_{+}} ( \sum_{r=1}^{m} |a_{r}|^{p}
)^{1/p} $ for \mbox{$\ell_{p}$-norm} in ${\mathbb{R}}^{m}$ it follows that
$|\theta_{\ell}| \leq2^{j_{0}(1/2-1/p)_{+}} ( \sum_{k
=0}^{2^{j_{0}}-1} |c_{j_{0},k}|^{p} )^{1/p} +\break \sum_{j =
j_{0}}^{+ \infty} 2^{j(1/2-1/p)_{+}} ( \sum_{k=0}^{2^{j}-1}
|\beta_{j,k}|^{p} )^{1/p}$.

Since $f \in B^{s}_{p,q}(A)$, one has that $ ( \sum
_{k=0}^{2^{j}-1} |\beta_{j,k}|^{p} )^{1/p} \leq A 2^{-j(s + 1/2-1/p)}$
and
$ ( \sum_{k =0}^{2^{j_{0}}-1} |c_{j_{0},k}|^{p} )^{1/p}
\leq A$
which implies that\vspace*{1pt}
$|\theta_{\ell}| \leq A 2^{j_{0}(1/2-1/p)_{+}} + A\times\break \sum_{j =
j_{0}}^{+ \infty} 2^{-j(s + 1/2-1/p - (1/2-1/p)_{+})}$.
Taking, for\vspace*{1pt} instance, $j_{0} = 0$ completes the proof since by assumption
$s + 1/2-1/p > 0$.
\end{pf}

\textit{An upper bound for ${\mathbb E}| \hat{\beta}_{j,k} -
\beta
_{j,k}|^{2}$} (the proof to control ${\mathbb E}|\hat{c}_{j_{0},k} -
c_{j_{0},k} |^{2} $ follows from the same arguments): from the
decomposition (\ref{eq:Z1Z2}) it follows that ${\mathbb E}| \hat
{\beta
}_{j,k} - \beta_{j,k}|^{2} \leq2 {\mathbb E}|Z_{1}|^{2} + 2 {\mathbb E}
|Z_{2}|^{2}$. Since $\eta_{\ell}$ are i.i.d. $N_{{\mathbb{C}}}
(0,1/n )$, the bound (\ref{eq:bound1}) implies that
%
%
\begin{equation} \label{eq:boundZ22}
{\mathbb E}|Z_{2}|^{2} \leq C \frac{ 2^{2 j \nu}}{n}.
\end{equation}
Then let us write $Z_{1} = \frac{1}{n} \sum_{m=1}^{n} (W_{m} -
{\mathbb E}
W_{m} )$ with $W_{m} = h_{j,k}( \tau_m)$ and $h_{j,k}(\tau) = \sum
_{\ell\in\Omega_{j}} \frac{\psi^{j,k}_{\ell} \theta_{\ell} }{
\gamma_{\ell}} e^{-i 2 \pi\ell\tau}$ for $\tau\in{\mathbb{R}}$. By
independence of the $\tau_{m}$'s, one has that ${\mathbb E}|Z_{1}|^{2}
\leq
\frac{1}{n} {\mathbb E}| W_{1}|^{2}$. Applying Lemma \ref{lem:pars}
with $h
= h_{j,k}$ and since the density $G$ is bounded, it follows that
%
%
\begin{eqnarray}\label{eq:W1}
{\mathbb E}| W_{1}|^{2} & = & \int_{{\mathbb{R}}} |h_{j,k}(\tau
)|^{2} g(\tau) \,d
\tau\nonumber\\
&=& \int_{0}^{1} |h_{j,k}(\tau)|^{2} G(\tau) \,d
\tau
\leq C \int_{0}^{1} |h_{j,k}(\tau)|^{2} \,d \tau\\
&\leq& C \sum
_{\ell\in\Omega_{j}} \frac{|\psi^{j,k}_{\ell}|^2 |\theta_{\ell
}|^2}{ |\gamma_{\ell}|^2 },\nonumber
\end{eqnarray}
where the last inequality follows from Parseval's relation. Then, using
the bound (\ref{eq:bound1}) and Lemma \ref{lem:bound}, inequality
(\ref{eq:W1}) implies that there exists a constant $C$ such that for
all $f \in B^{s}_{p,q}(A)$
%
%
\begin{equation} \label{eq:boundZ12}
{\mathbb E}Z_{1}^{2} \leq C \frac{1}{n} \sum_{\ell\in\Omega_{j}}
\frac
{|\psi^{j,k}_{\ell}|^2 }{ |\gamma_{\ell}|^2 } \leq C \frac{2^{2j
\nu}}{n}.
\end{equation}
Hence, using the bounds (\ref{eq:boundZ22}) and (\ref{eq:boundZ12}),
it follows that there exists a constant $C$ such that for all $f \in
B^{s}_{p,q}(A)$,
${\mathbb E}| \hat{\beta}_{j,k} - \beta_{j,k}|^{2} \leq C \frac
{2^{2j \nu}}{n}$.

\textit{An upper bound for ${\mathbb E}| \hat{\beta}_{j,k} -
\beta
_{j,k}|^{4}$}: from the decomposition (\ref{eq:Z1Z2}) it follows that
${\mathbb E}| \hat{\beta}_{j,k} - \beta_{j,k}|^{4} \leq C ({\mathbb
E}|Z_{1}|^{4}
+ {\mathbb E}|Z_{2}|^{4})$. As $Z_{2}$ is a centered Gaussian variable with
variance $\frac{1}{n} \varepsilon^{2} \sum_{\ell\in\Omega_{j}} |
\frac{\psi^{j,k}_{\ell}}{ \gamma_{\ell}} |^{2} \leq C \frac{ 2^{2
j \nu}}{n}$, one has that
%
%
\begin{equation} \label{eq:boundZ24}
{\mathbb E}|Z_{2}|^{4} \leq C \frac{ 2^{ j 4 \nu}}{n^{2}}.
\end{equation}
Then remark that $Z_{1} = \frac{1}{n} \sum_{m=1}^{n}Y_{m}$ with
$Y_{m} = \sum_{\ell\in\Omega_{j}} \frac{\psi^{j,k}_{\ell} \theta
_{\ell} }{ \gamma_{\ell}} (e^{-i 2 \pi\ell\tau_{m}} - \gamma
_{\ell})$, and recall the so-called Rosenthal inequality for moment
bounds of i.i.d. variables \cite{MR0440354}: if $X_{1},\ldots,X_{n}$
are i.i.d. random variables such that ${\mathbb E}X_{j} = 0$, ${\mathbb
E}X_{j}^{2}
\leq\sigma^{2}$, there exists a positive constant $C$ such that
${\mathbb E}| \sum_{j=1}^{n} X_{j}/n |^{4} \leq C ( \sigma
^{4}/n^{2} +
{\mathbb E}|X_{1}|^{4}/n^{3})$.

Now remark that ${\mathbb E}Y_{m} = 0 $, and arguing as previously for the
control of ${\mathbb E}|W_{1}|^{2}$, see (\ref{eq:W1}), it follows
that ${\mathbb E}|Y_{m}|^{2} \leq C 2^{2 j\nu}$ where $C$ is constant not
depending on~$f$. Then remark that
\begin{eqnarray}
{\mathbb E}|Y_{1}|^{4} \leq C \biggl( \int_{{\mathbb{R}}}
|h_{j,k}(\tau)|^{4} g(\tau
) \,d \tau+ |\beta_{jk}|^{4} \biggr) \nonumber\\
\eqntext{\mbox{with } \displaystyle h_{j,k}(\tau) =
\sum_{\ell\in\Omega_{j}} \frac{\psi^{j,k}_{\ell} \theta_{\ell}
}{ \gamma_{\ell}} e^{-i 2 \pi\ell\tau},}
\end{eqnarray}
and that
\[
\int_{{\mathbb{R}}} |h_{j,k}(\tau)|^{4} g(\tau) \,d \tau\leq\sup
_{\tau\in
{\mathbb{R}}} \{ |h_{j,k}(\tau)|^{2} \} \int_{{\mathbb{R}}}
|h_{j,k}(\tau)|^{2}
g(\tau) \,d \tau.
\]
Note that\vspace*{-1pt} using (\ref{eq:bound1}) and Lemma \ref{lem:bound}, it
follows that $|h_{j,k}(\tau)| \leq\sum_{\ell\in\Omega_{j}} \frac
{|\psi^{j,k}_{\ell}| | \theta_{\ell} | }{| \gamma_{\ell}|}\leq
C \sum_{\ell\in\Omega_{j}} \frac{|\psi^{j,k}_{\ell}| }{| \gamma
_{\ell}|} \leq C 2^{j (\nu+1/2)} $ uniformly\vspace*{1pt} for $f \in
B^{s}_{p,q}(A)$. Then, arguing as for the control of ${\mathbb E}|
W_{1}|^{2}$, see (\ref{eq:W1}), one has that
$\int_{{\mathbb{R}}} |h_{j,k}(\tau)|^{2} g(\tau) \,d \tau\leq C 2^{2j
\nu}$, which finally implies that
${\mathbb E}|Y_{1}|^{4} \leq C 2^{j(4 \nu+ 1)}$,
since $ |\beta_{jk}| \leq C$ uniformly for $f \in B^{s}_{p,q}(A)$.
Then, using Rosenthal's inequality, it follows that there exists a
constant $C$ such that for all $f \in B^{s}_{p,q}(A)$
%
%
\begin{equation}\label{eq:boundZ14}
{\mathbb E}|Z_{1}|^{4} \leq C \biggl(\frac{ 2^{j 4 \nu}}{n^{2}} + \frac{
2^{j (4
\nu+ 1)}}{n^{3}}\biggr),
\end{equation}
which completes the proof for the control of ${\mathbb E}| \hat{\beta}_{j,k}
- \beta_{j,k}|^{4}$ using (\ref{eq:boundZ24}) and (\ref
{eq:boundZ14}).
\end{pf*}


\subsection{\texorpdfstring{Proof of Proposition \protect\ref{prop:assproba} and
Theorem \protect\ref{theo:upperbound}}{Proof of Proposition 4 and Theorem 2}}

First, let us prove the following proposition.
%
%
\begin{prop} \label{prop:proba} Let $f \in L^{2}([0,1])$, $n \geq1$
and $j \geq0$. Suppose that $g$ satisfies Assumption \ref{assdecay}.
For $0 \leq k \leq2^{j}-1$ and $\theta_{\ell} = \int_{0}^{1} f(x)
e^{-i2\pi\ell x}\,dx$, define
\[
\sigma^{2}_{j} = {2^{-j} \varepsilon^{2} \sum_{\ell\in\Omega_{j}}} |
\gamma_{\ell} |^{-2},\qquad
V^{2}_{j} = {\| g \|_{\infty} 2^{-j} \sum_{\ell\in\Omega_{j}}} \frac
{|\theta_{\ell}|^{2}}{ |\gamma_{\ell}|^{2}}
\]
and
\[
\delta_{j} = {2^{-j/2} \sum_{\ell\in\Omega_{j}}} \frac{|\theta
_{\ell}|}{ |\gamma_{\ell}| }
\]
with $\| g \|_{\infty} = \sup_{x \in{\mathbb{R}}} \{ g(x) \}$. Let
$t >0$, then,
\[
{\mathbb P} \Biggl( |\hat{\beta}_{j,k} - \beta_{j,k} | \geq2 \max
\Biggl(
\sigma_{j} \sqrt{\frac{2 t}{n}} , \sqrt{\frac{ 2 V^{2}_{j}t}{n}
} + \delta_{j} \frac{t}{3n} \Biggr) \Biggr) \leq2 \exp(-t).
\]
\end{prop}
\begin{pf}
Let $u > 0$, and remark that from the decomposition (\ref{eq:Z1Z2}) it follows
\[
{\mathbb P} ( |\hat{\beta}_{j,k} - \beta_{j,k} | \geq u
) \leq
{\mathbb P}( |Z_{1} | \geq u/2) + {\mathbb P}( |Z_{2} | \geq u/2).
\]
Recall that the $\eta_{\ell}$'s are i.i.d. $N_{{\mathbb{C}}}
(0,1/n )$. Hence, $Z_{2}$ is a centered Gaussian variable with
variance $\frac{1}{n} \varepsilon^{2} \sum_{\ell\in\Omega_{j}} |
\frac{\psi^{j,k}_{\ell}}{ \gamma_{\ell}} |^{2} \leq\frac{1}{n}
\sigma^{2}_{j}$, with $\sigma^{2}_{j} = 2^{-j} \varepsilon^{2} \sum
_{\ell\in\Omega_{j}} | \gamma_{\ell} |^{-2}$, which implies that
(see, e.g., \cite{massart}) for any $t > 0$
%
%
\begin{equation}\label{eq:Z2}
{\mathbb P}\Biggl( |Z_{2}| \geq\sigma_{j} \sqrt{\frac{2 t}{n}} \Biggr) \leq2
\exp(-t ).
\end{equation}
By definition,
\[
\tilde{\gamma}_{\ell} = \frac{1}{n} \sum_{m=1}^n
e^{-i 2 \pi\ell\tau_m},
\]
and thus $Z_{1} = \frac{1}{n} \sum
_{m=1}^{n} (W_{m} - {\mathbb E}W_{m} )$ with $W_{m} = {\sum_{\ell\in
\Omega
_{j}} }\frac{\psi^{j,k}_{\ell} \theta_{\ell} }{ \gamma_{\ell}}
e^{-i 2 \pi\ell\tau_m}$. Remark that $W_{m}$ are random variables
bounded by $\delta_{j} = {2^{-j/2} \sum_{\ell\in\Omega_{j}} }\frac{
|\theta_{\ell}|}{ |\gamma_{\ell}| } $. Moreover, using Lemma \ref
{lem:pars} with $h = h_{j,k}(\tau) = \sum_{\ell\in\Omega_{j}}
\frac{\psi^{j,k}_{\ell} \theta_{\ell} }{ \gamma_{\ell}} e^{-i 2
\pi\ell\tau}$ for $\tau\in{\mathbb{R}}$ it follows that
\[
{\mathbb E}|W_{1}|^{2} = \int_{{\mathbb{R}}} |h_{j,k}(\tau)|^{2}
g(\tau) \,d \tau
\leq\| G \|_{\infty} \sum_{\ell\in C_{j}} \frac{|\psi^{j,k}_{\ell
}|^2 |\theta_{\ell}|^2}{ |\gamma_{\ell}|^2 } \leq V_{j}^{2},
\]
where $ V_{j}^{2} = \| g \|_{\infty} 2^{-j} \sum_{\ell\in C_{j}}
\frac{ |\theta_{\ell}|^2}{ |\gamma_{\ell}|^2 }$, since $|\psi
^{j,k}_{\ell}| \leq2^{-j/2}$ and $ \| g \|_{\infty} = \| G \|
_{\infty}$. Hence, from Bernstein's inequality it follows that for any
$t > 0$ (see, e.g., Proposition 2.9 in \cite{massart})
%
%
\begin{equation} \label{eq:Z1}
{\mathbb P} \Biggl( |Z_{1}| \geq\sqrt{\frac{ 2 V_{j}^{2} t }{n} } +
\delta
_{j} \frac{t }{3n} \Biggr) \leq2 \exp(-t).
\end{equation}
Taking $u = 2 \max( \sigma_{j} \sqrt{\frac{2 t}{n}} ,
\sqrt{\frac{ 2 V_{j}^{2} t }{n} } + \delta_{j} \frac{t }{3n}
)$ for $t > 0$ concludes the proof of Proposition \ref{prop:proba}.
\end{pf}


Proposition \ref{prop:proba} would suggest to take level-dependent
thresholds of the form
%
%
\begin{equation} \label{eq:thr}\quad
\lambda_{j,k}^{\ast} = \lambda_{j}^{\ast} = 2 \max\Biggl( \sigma
_{j} \sqrt{\frac{2 \eta\log(n) }{n}} , \sqrt{\frac{ 2 \eta
V^{2}_{j} \log(n) }{n} } + \delta_{j} \frac{ \eta\log(n) }{3n} \Biggr)
\end{equation}
for some constant $\eta> 0$. The first term in the maximum
(\ref{eq:thr}) is the classical universal threshold with
heteroscedastic\vspace*{-1pt}
variance $ \sigma_{j}^{2}$ which corresponds to an upper bound of the
variance of the Gaussian term $\varepsilon\sum_{\ell\in\Omega_{j}}
\frac{\eta_{l}}{\gamma_{\ell}} $ in the expression of~$\hat{\beta
}_{j,k}$. However, the second term in the maximum (\ref{eq:thr})
depends on the modulus of the unknown Fourier coefficients $\theta
_{\ell}$, and thus the thresholds $\lambda_{j}^{\ast}$ cannot be
used in practice.

Fortunately, the computation of the threshold $\lambda_{j,k}^{\ast}$
can be simplified using the following arguments. Since there exist two
constants $C_{1},C_{2}$ such that for all $\ell\in\Omega_{j}$,
$C_{1} 2^{j} \leq\ell\leq C_{2} 2^{j}$ and since $\lim_{|\ell| \to
+ \infty} \theta_{\ell} = 0$ uniformly for $f \in B_{p,q}^{s}(A)$
it follows that as $j \to+ \infty$
\[
V^{2}_{j} = \| g \|_{\infty} 2^{-j} \sum_{\ell\in\Omega_{j}} \frac
{|\theta_{\ell}|^{2}}{ |\gamma_{\ell}|^{2}} = o \biggl( 2^{-j} \sum
_{\ell\in\Omega_{j}} |\gamma_{\ell}|^{-2} \biggr) = o (
\sigma^{2}_{j} ).
\]
Also, if $f \in B_{p,q}^{s}(A)$ then $\sum_{\ell\in\Omega_{j}}
|\theta_{\ell}|^{2} = o(1)$ as $j \to+ \infty$, and thus by the
Cauchy--Schwarz inequality, then as $j \to+ \infty$,
\[
\delta_{j} = 2^{-j/2} \sum_{\ell\in\Omega_{j}} \frac{|\theta
_{\ell}|}{ |\gamma_{\ell}| } \leq2^{-j/2} \biggl( \sum_{\ell\in
\Omega_{j}} |\theta_{\ell}|^{2} \biggr)^{1/2} \biggl( \sum_{\ell
\in\Omega_{j}} |\gamma_{\ell}|^{-2} \biggr)^{1/2} = o (
\sigma_{j} ),
\]
which finally implies that $ \lambda_{j}^{\ast} = o ( \lambda
_{j} ) $ as $j \to+ \infty$ where
$\lambda_{j} = 2 \sigma_{j} \sqrt{\frac{2 \eta\log(n)}{n}}$.
Hence, if one chooses $j_{0}$ to be slowly growing with $n$ [e.g.,
$j_{0} = \log(\log(n))$], or avoid thresholding at very low
resolution levels, then the threshold $\lambda_{j}$ can be used
instead of $ \lambda_{j}^{\ast}$ whose computation would require an
estimator of the $|\theta_{\ell}|$'s.

Combining Propositions \ref{prop:moments} and \ref{prop:proba}, the
above remarks on the thresholds $\lambda_{j}$ and~$ \lambda_{j}^{\ast
}$, and by arguing as in the proof of Theorem 1 in \cite{JKPR04jrssb},
then Proposition \ref{prop:assproba} and Theorem \ref
{theo:upperbound} follow.


\subsection{\texorpdfstring{Proof of Theorem
\protect\ref{theo:lowerbound}}{Proof of Theorem 3}}

Let us fix a resolution $j \geq0$ whose choice will be discussed later
on, and consider for any $\eta=(\eta_i)_{i \in\{0,
\ldots,2^j-1\}} \in\{\pm1\}^{2^j}$ the function $f_{j,\eta}$ defined as
$f_{j,\eta} = \gamma_j \sum_{i=0}^{2^j-1} \eta_k \psi_{j,k}$,
where $\gamma_j = c 2^{-j(s+1/2)}$,
and $c$ is a positive constant satisfying $c \leq A$ which implies that
$f_{j,\eta} \in
B_{p,q}^s(A)$. For some $0 \leq i \leq2^{j} -1$ and $\eta\in\{\pm1\}
^{2^j}$, define also the vector $\eta^i \in\{\pm1\}^{2^j}$ with
components equal to those of $\eta$ except the $i$th one.

Let $\psi_{j,k} \star g(x) = \int_{{\mathbb{R}}} \psi_{j,k}(x-u)
g(u) \,du$.
By Parseval's relation, one has that $\|\psi_{j,k} \star g \|^{2} =
{\sum_{\ell\in\Omega_{j} }} |\psi^{j,k}_{\ell}|^{2} |\gamma_{\ell
}|^{2}$. Hence, under Assumption \ref{assordi} of a polynomial decay
for $\gamma_{\ell}$ and using the fact that $|\psi^{j,k}| \leq
2^{-j/2}$ for Meyer wavelets (see \cite{JKPR04jrssb}) it follows that
there exists a constant $C$ such that
$\|\psi_{j,k} \star g \|^{2} \leq C 2^{-2 j \nu}$.

\subsubsection{Asymptotic settings} \label{sec:alsetting-reg}

We set the resolution\vspace*{2pt} $j = j(n)$ to be the largest integer satisfying
$2^{j(n)} \leq n^{{1}/({2s + 2 \nu+1})}$.
However, to simplify the presentation, the dependency of $j$ on $n$ is
dropped in what follows. The definition of $f_{j,\eta}$, $\gamma_j $
and the bound $\|\psi_{j,k} \star g \|^{2} \leq C 2^{-2 j \nu}$ thus
imply that
\[
\gamma_j  = \mathcal O\bigl( n^{- ({s+1/2})/({2s + 2 \nu+1}) }\bigr)
\]
and
\begin{eqnarray*}
\|f_{j,\eta}\|^2 &=& \mathcal O\bigl( n^{- ({2s})/({2s+2\nu+1})} \bigr),
\\
\|f_{j,\eta} \star g\|^2 & = & \biggl\| \gamma_j \sum_k \eta_k (\psi
_{j,k} \star g)
\biggr\|^2 =\mathcal O \bigl( n^{({-2s-2\nu})/({2s+2\nu+1})} \bigr), \\
\|(f_{j,\eta} - f_{j,\eta^i}) \star g\|^2 & = & \|2 \gamma_j \eta_i
(\psi_{j,i}
\star g) \|^2 =\mathcal O ( \gamma_j^2 2^{-2 j \nu} ) =
\mathcal O
(n^{-1} ).
\end{eqnarray*}
From the above equations, we can thus conclude that
$n \|(f_{j,\eta} - f_{j,\eta^i}) \star g\|^2 = \mathcal O (
1 )$,
but note that the term $n \|f_{j,\eta} \star g\|^2$ does not converge
to $0$. At last, observe that by assumption $s > 2\nu+1$ which implies
that $n \|f_{j,\eta} \star g\|^3 \to0$, $n \|(f_{j,\eta} - f_{j,\eta
^i}) \star g\| \|f_{j,\eta}\| \|f_{j,\eta}\star g\| \to0$ and $n \|
f_{j,\eta}\|^3 \to0$.

\subsubsection{Likelihood ratio}

Let $F(Y)$ be real valued and bounded measurable function of the $n$
trajectories
$Y=(Y_{1},\ldots,Y_{n})$. Because of the independence of the $\tau
_i$'s and the $W_i$'s, we have that
\begin{eqnarray*}
\mathbb{E}_{f} [ F(Y) ] & =& \int_{{\mathbb{R}}^{n}}
\mathbb{E}_{f,W}
[F(Y) | \tau_1=t_{1},\ldots,\tau_{n}=t_{n} ] g(t_{1}) \,d
t_{1} \cdots g(t_{n}) \,d t_{n},
\end{eqnarray*}
where $\mathbb{E}_{f}$ denotes the expectation with respect to the law
of $Y=(Y_{1},\ldots,Y_{n})$
when $f$ is the true hypothesis, and $\mathbb{E}_{f,W}$ is used to
denote expectation only with
respect to law of the Brownian motions $W_{1},\ldots,W_{n}$ where the
shifts are fixed and $f$ is the
true hypothesis. Now using the classical Girsanov formula it follows
that for any function $h \in L^{2}([0,1])$
\begin{eqnarray*}
\mathbb{E}_{f} [ F(Y) ] & =& \int_{{\mathbb{R}}^{n}}
\mathbb{E}_{h,W}
[F(Y) | \tau_1=t_{1},\ldots,\tau_{n}=t_{n} ] \Lambda
_{n}(f,h) g(t_{1}) \,d t_{1}\cdots g(t_{n}) \,d t_{n} \\
& = & \mathbb{E}_{h} [F(Y) \Lambda_{n}(f,h) ],
\end{eqnarray*}
where $\Lambda_{n}(f,h)$ is the following likelihood ratio:
\[
\Lambda_{n}(f,h) = \prod_{i=1}^{n} \exp\biggl( \int_{0}^{1}
\bigl(f(x-\tau_{i})-h(x-\tau_{i})\bigr) \, d Y_{i}(x) + \frac{1}{2} \|h\|^{2} -
\frac{1}{2}
\|f\|^{2} \biggr).
\]
In what follows, $f_{0}$ is used to denote the hypothesis $f \equiv0$.

\subsubsection{Technical lemmas}

Given $n$ arbitrary trajectories $Y_1, \ldots, Y_n$ from model (\ref
{model}), we define $\mathbb{E}_{\tau}
( \Lambda_{n}(f_{j, \eta^{i}},f_{0}) )$ as the
expectation of the likelihood ratio with respect to the law of the
random shifts, namely
\begin{eqnarray*}
{\mathbb E}_{\tau} ( \Lambda_{n}(f_{j, \eta^{i}},f_{0})
) & = &
\int_{\mathbb{R}^n} \prod_{i=1}^{n} e^{ \int_{0}^{1}
(f(x-\tau_{i})-h(x-\tau_{i}))\, d Y_{i}(x) + {1}/{2} \|h\|^{2} -
{1}/{2}
\|f\|^{2}} \\
&&\hspace*{29.7pt}{}\times g(\tau_1) \cdots g(\tau_n)\, d\tau_1 \cdots d\tau_n.
\end{eqnarray*}
\begin{lemma}\label{lem:assouad-bis} Suppose for some constants
$\lambda> 0$ and $\pi_0 > 0 $ and all sufficiently large $n$ we have that
%
%
\begin{equation}\label{eq:min-ratio}
{\mathbb P}_{f_{j, \eta}} \biggl( \frac{ {\mathbb E}_{\tau} (
\Lambda_{n}(f_{j,
\eta^{i}},f_{0}) )}
{ {\mathbb E}_{\tau} ( \Lambda_{n}(f_{j, \eta},f_{0})
)} \geq
e^{-\lambda} \biggr) \geq
\pi_{0}
\end{equation}
for all $f_{j,\eta}$ and all $i \in\{0 ,\ldots,2^{j}-1\}$. Then there
exists a
positive constant $C$, such that for all sufficiently large $n$ and any
estimator $\hat{f}_n$
\[
\max_{\eta\in\{\pm1\}^{2^j}} {\mathbb E}_{f_{j, \eta}} \|\hat
{f}_{n} -
f_{j,\eta}
\|^{2} \geq C \pi_0 e^{- \lambda} 2^{j} \gamma_{j}^{2}.
\]
\end{lemma}
\begin{pf}
Our proof is inspired by the proof of Lemma 2.10 in \cite{HKPT}. For
this, let $I_{jk} = [\frac{k}{2^j},\frac{k+1}{2^j}]$ and arguing as in
\cite{HKPT} it follows that for any estimator $\hat{f}_{n}$
\begin{eqnarray*}
&&\max_{\eta\in\{\pm1\}^{2^j}} {\mathbb E}_{f_{j, \eta}} \|\hat
{f}_{n} -
f_{j,\eta} \|^{2} \\
&&\qquad \geq 2^{-2^{j}}
\sum_{k = 0}^{2^{j}-1} \sum_{\eta| \eta_{k} = 1} {\mathbb E}_{f_{j,
\eta
}} \biggl[ \int_{I_{j,k}}
|\hat{f}_{n} - f_{j,\eta}|^{2} \\
& &\qquad\quad\hspace*{95.48pt}{} + \Lambda_{n}(f_{j, \eta^{k}},f_{j, \eta}) \int
_{I_{j,k}} |\hat{f}_{n} - f_{j,\eta^{k}}|^{2} \biggr].
\end{eqnarray*}
Let $Z(Y) = [ \int_{I_{j,k}} |\hat{f}_{n} - f_{j,\eta}|^{2} +
\Lambda_{n}(f_{j, \eta^{k}},f_{j, \eta})
\int_{I_{j,k}} |\hat{f}_{n} - f_{j,\eta^{k}}|^{2} ]$ and
remark that
\begin{eqnarray*}
{\mathbb E}_{f_{j, \eta}} [ Z(Y) ] & = & {\mathbb E}_{f_{0},W} \int
_{{\mathbb{R}}^{n}}
\biggl[ \Lambda_{n}(f_{j, \eta},f_{0}) \int_{I_{j,k}} |\hat{f}_{n}
- f_{j,\eta}|^{2} \\
& &\hspace*{48.2pt}{} + \Lambda_{n}(f_{j, \eta^{k}},f_{0}) \int_{I_{j,k}}
|\hat{f}_{n} - f_{j,\eta^{k}}|^{2} \biggr]\\
& &\hspace*{42.5pt}{}\times  g(\tau_{1}) \,d \tau
_{1}\cdots g(\tau_{n}) \,d \tau_{n}.
\end{eqnarray*}
Now, since under the hypothesis $f_{0}$ the trajectories $Y_{1},\ldots,Y_{n}$
do not depend on the random shifts
$\tau_{1},\ldots,\tau_{n}$ it follows that
$\hat{f}_{n}$ does not depend on the shifts $\tau_{1},\ldots,\tau
_{n}$ as it
is by definition a measurable function with respect to the sigma algebra
generated by $Y_{1},\ldots,Y_{n}$. This implies that for any $\delta> 0$
\begin{eqnarray*}
{\mathbb E}_{f_{j, \eta}} [ Z(Y) ] & = & {\mathbb E}_{f_{0},W}
\biggl[{\mathbb E}_{\tau}
( \Lambda_{n}(f_{j, \eta},f_{0}) )
\int_{I_{j,k}} |\hat{f}_{n} - f_{j,\eta}|^{2} \\
& &\hspace*{30.66pt} {} + {\mathbb E}_{\tau} ( \Lambda_{n}(f_{j, \eta^{k}},f_{0})
)
\int_{I_{j,k}} |\hat{f}_{n} - f_{j,\eta^{k}}|^{2} \biggr] \\
& \geq& {\mathbb E}_{f_{0},W} \bigl[{\mathbb E}_{\tau} (
\Lambda_{n}(f_{j,
\eta},f_{0}) ) \delta^{2} \mathbh{1}_{ \{ \int_{I_{j,k}}
|\hat{f}_{n} - f_{j,\eta}|^{2} > \delta^{2} \} } \\
& &\hspace*{28.8pt} {}+ {\mathbb E}_{\tau} ( \Lambda_{n}(f_{j, \eta^{k}},f_{0}))
\delta^{2} \mathbh{1}_{ \{ \int_{I_{j,k}} |\hat{f}_{n} -
f_{j,\eta
^{k}}|^{2} > \delta^{2} \} }
\bigr].
\end{eqnarray*}
Now, remark that
\begin{eqnarray*}
&&\biggl( \int_{I_{j,k}} |\hat{f}_{n} - f_{j,\eta}|^{2} \biggr)^{1/2}
+ \biggl(\int_{I_{j,k}} |\hat{f}_{n} - f_{j,\eta^{k}}|^{2}
\biggr)^{1/2} \\
&&\qquad \geq \biggl(\int_{I_{j,k}} | f_{j,\eta} - f_{j,\eta
^{k}}|^{2} \biggr)^{1/2}
\\
&&\qquad\geq 2 \gamma_{j} \biggl(\int_{I_{j,k}} |\psi_{jk} |^{2} \biggr)^{1/2},
\end{eqnarray*}
and let us argue as in the proof of Lemma 2 in \cite{willer} to a find
a lower bound for $\int_{I_{j,k}} |\psi_{jk} |^{2}$. By definition
(see Section \ref{sec:estimfnonlinear}) $ \psi_{j,k}(x) = 2^{j/2} \sum
_{i \in{\mathbb{Z}}} \psi^{\ast}(2^{j}(x+i)-k)$ where $ \psi^{\ast
}$ is
the Meyer wavelet over ${\mathbb{R}}$ used to construct $\psi$. A
change of
variable shows that $\int_{I_{j,k}} |\psi_{jk}(x) |^{2} \,dx = \int
_{0}^{1} | \sum_{i \in{\mathbb{Z}}} \psi^{\ast}(x+2^{j} i)
|^{2} \,dx$ which implies that $\int_{I_{j,k}} |\psi_{jk}(x) |^{2} \,dx
\geq\int_{0}^{1} | \psi^{\ast}(x) |^{2} \,dx - \sum_{i
\in{\mathbb{Z}}^{\ast}} \int_{0}^{1} | \psi^{\ast}(x+2^{j}
i)
|^{2} \,dx $. Now as $\psi^{\ast}$ has a fast decay, it follows that
there exists a constant $A > 0$ such that $ | \psi^{\ast}(x)| \leq
\frac{A}{1+x^{2}}$. Thus, $\int_{I_{j,k}} |\psi_{jk}(x) |^{2} \,dx
\geq\int_{0}^{1} | \psi^{\ast}(x) |^{2} \,dx - A^{2}
2^{-2j} \sum_{i \in{\mathbb{Z}}^{\ast}} i^{-2}$. Hence, it follows that
there exists a constant $\rho> 0$ such that $ (\int_{I_{j,k}}
|\psi_{jk} |^{2} )^{1/2} \geq\rho$ for any $k$, and all $j$
sufficiently large.

Hence, if one takes $\delta= 2 \rho\gamma_{j}$ it follows that
\[
\mathbh{1}_{ \{ \int_{I_{j,k}} |\hat{f}_{n} - f_{j,\eta}|^{2}
> \delta
^{2} \} } \geq\mathbh{1}_{ \{ \int_{I_{j,k}} |\hat
{f}_{n} -
f_{j,\eta^{k}}|^{2} \leq\delta^{2} \} },
\]
which yields
\begin{eqnarray*}
{\mathbb E}_{f_{j, \eta}} [ Z(Y) ] & \geq& \delta^{2} {\mathbb
E}_{f_{0},W}
\biggl[ {\mathbb E}_{\tau} ( \Lambda_{n}(f_{j, \eta},f_{0}) )
\min\biggl(1, \frac{ {\mathbb E}_{\tau} ( \Lambda_{n}(f_{j,
\eta
^{k}},f_{0}) )} { {\mathbb E}_{\tau}
( \Lambda_{n}(f_{j, \eta},f_{0}) )} \biggr) \biggr] \\
& = & \delta^{2} {\mathbb E}_{f_{0}} \biggl[ \Lambda_{n}(f_{j, \eta},f_{0})
\min\biggl(1, \frac{ {\mathbb E}_{\tau} ( \Lambda_{n}(f_{j,
\eta
^{k}},f_{0}) )} { {\mathbb E}_{\tau}
( \Lambda_{n}(f_{j, \eta},f_{0}) )} \biggr) \biggr] \\
& = & \delta^{2} {\mathbb E}_{f_{j, \eta}} \biggl[ \min\biggl(1,
\frac{
{\mathbb E}_{\tau}
( \Lambda_{n}(f_{j, \eta^{k}},f_{0}) )} { {\mathbb
E}_{\tau}
( \Lambda_{n}(f_{j, \eta},f_{0}) )} \biggr) \biggr],
\end{eqnarray*}
and arguing as in the proof of Lemma 2.10 in \cite{HKPT} completes the
proof.
\end{pf}


Now remark that under the hypothesis $f = f_{j, \eta}$, then
each $Y_i$ is
given by $d Y_i(x) = f_{j,\eta}(x-\alpha_i) \,dx + dW_i(x)$ where each
$\alpha_i$ is the true random shift of the $i$th trajectory. Thus, under
this hypothesis, we obtain
\begin{eqnarray*}
&&{\mathbb E}_{\tau} ( \Lambda_{n}(f_{j, \eta}, f_{0})
)\\
&&\qquad=\prod_{i=1}^{n}
\int_{{\mathbb{R}}} g(\tau_{i}) e^{ [ \int_{0}^{1} f_{j, \eta
}(x-\tau
_{i})f_{j,
\eta}(x-\alpha_{i}) \,dx + f_{j, \eta}(x-\tau_{i}) \,d W_{i}(x) - {1}/{2}
\|f_{j, \eta}\|^{2} ]} \,d \tau_{i}
\end{eqnarray*}
and
\begin{eqnarray*}
&&{\mathbb E}_{\tau} ( \Lambda_{n}(f_{j, \eta^{i}},f_{0}))\\
&&\qquad =
\prod_{i=1}^{n} \int_{{\mathbb{R}}} g(\tau_{i}) e^{ [ \int
_{0}^{1} f_{j,
\eta^{i}}(x-\tau_{i}) f_{j, \eta}(x-\alpha_{i}) \,dx + f_{j,
\eta^{i}}(x-\tau_{i}) \,d W_{i}(x) - {1}/{2} \|f_{j, \eta^{i}}\|
^{2} ]} \,d \tau_{i}.
\end{eqnarray*}
Using the two expressions above, we now study condition (\ref
{eq:min-ratio}).
\begin{lemma}\label{lem:min-proba}
Following the choices of $j(n)$ and $\gamma_{j(n)}$ given in our
algebraic setting, there exist $\lambda>0$
and $\pi_0>0$ such that for all sufficiently large $n$
\[
{\mathbb P}_{f_{j, \eta}} \biggl( \frac{ {\mathbb E}_{\tau} (
\Lambda_{n}(f_{j,
\eta^{i}},f_{0}) )}
{ {\mathbb E}_{\tau} ( \Lambda_{n}(f_{j, \eta},f_{0})
)} \geq
e^{-\lambda} \biggr) \geq\pi_{0}.
\]
\end{lemma}
\begin{pf}
To obtain the required bound, we use several second order Taylor
expansions. From the Cauchy--Schwarz inequality, we have
\begin{eqnarray*}
e^{ \int_{0}^{1} f_{j, \eta}(x-\tau_{i}) f_{j,
\eta}(x-\alpha_{i}) \,dx} &=& 1+ \int_{0}^{1} f_{j, \eta}(x-\tau_{i}) f_{j,
\eta}(x-\alpha_{i}) \,dx \\
&&{}+ \mathcal O_p(\|f_{j,\eta}\|^4).
\end{eqnarray*}
A similar argument yields
$\mathbb{E} [ | \int_{0}^1 f_{j, \eta}(x-\tau_i) \,dW_i(x)
|] \leq\|f_{j, \eta}\|$,
and the Markov inequality used with a second order expansion implies
$e^{ \int_{0}^1 f_{j, \eta}(x-\tau_i) \,dW_i(x) } = 1+
\int_{0}^1
f_{j, \eta}(x-\tau_i) \,dW_i(x) + \frac{1}{2} [\int_{0}^1 f_{j,
\eta}(x-\tau_i) \,dW_i(x) ]^2 + \mathcal O_p(\|f_{j,\eta}\|^3)$.
Looking now at the complete expression of ${\mathbb E}_{\tau} (
\Lambda
_{n}(f_{j,
\eta},f_{0}) )$, we obtain
${\mathbb E}_{\tau} ( \Lambda_{n}(f_{j, \eta},f_{0}) ) =
\prod_{i=1}^n
\int_{\mathbb{R}} g(\tau_i) [1+ \int_{0}^{1} f_{j, \eta
}(x-\tau_{i}) f_{j,
\eta}(x-\alpha_{i}) \,dx +\int_{0}^1 f_{j, \eta}(x-\tau_i) \,dW_i(x)
+\break \frac{1}{2} [\int_{0}^1 f_{j, \eta
}(x-\tau_i)
\,dW_i(x) ]^2 - \frac{1}{2} \|f_{j,\eta}\|^2 + \mathcal O_p(\|
f_{j,\eta}\|^3)
]$.
The Fubini-type\vspace*{1pt} theorem for stochastic integrals (see,
e.g.,
\cite{IkeWat}, Chapter 3, Lemma 4.1) enables us to write
$
\log{\mathbb E}_{\tau} ( \Lambda_{n}(f_{j, \eta},f_{0})
) =
\sum_{i=1}^n \log[ 1+ \int_{0}^1 (f_{j,\eta}\star g)(x)
f_{j,\eta}(x-\alpha_i) \,dx + \int_{0}^1
(f_{j,\eta}\star g)(x) \,dW_i(x) \frac{1}{2} \int_{\mathbb{R}}
[\int_{0}^1 f_{j,
\eta}(x-\tau_i) \,dW_i(x) ]^2 g(\tau_i) \,d \tau_i - \frac{1}{2}
\|f_{j,\eta}\|^2 +\break \mathcal O_p(\|f_{j,\eta}\|^3)
]$.

Then applying a classical expansion of the logarithm $\log(1+z) = z -
\frac{z^2}{2} + \mathcal O(z^3)$, we obtain
%
%
\begin{eqnarray}\qquad
&&\log{\mathbb E}_{\tau}\Lambda_n(f_{j, \eta},f_0)\nonumber\\
&&\qquad = z - \frac
{z^2}{2} +
\mathcal O_p(z^3)
\nonumber\\
\label{eq:line1}
&&\qquad = \sum_{i=1}^n \int_{0}^1 (f_{j, \eta}\star g)(x)f_{j,\eta
}(x-\alpha_i) \,dx +
\int_{0}^1 (f_{j, \eta}\star g)(x) \,dW_i(x)\\
\label{eq:line2}
&&\qquad\quad{} + \frac{1}{2} \sum_{i=1}^n
\int_{\mathbb{R}} g(\tau_i) \biggl[\int_{0}^1
f_{j,\eta}(x-\tau_i) \,dW_i(x) \biggr]^2 - \frac{n}{2} \|f_{j,\eta}\|
^2\\[-0.5pt]
\label{eq:line3}
&&\qquad\quad{} - \frac{1}{2} \sum_{i=1}^n \biggl( \int_{0}^1 (f_{j,
\eta}\star g)(x)f_{j,\eta}(x-\alpha_i) \,dx\nonumber\\[-9pt]\\[-9pt]
&&\qquad\quad\hspace*{41.1pt}{} + \int_{0}^1 (f_{j, \eta
}\star g)(x)
\,dW_i(x) \nonumber\\[-0.5pt]
\label{eq:line4}
&&\qquad\quad\hspace*{41.1pt}{} +\frac{1}{2} \int_{\mathbb{R}} g(\tau_i) \biggl[\int
_{0}^1 f_{j,\eta}(x-\tau_i)
\,dW_i(x) \biggr]^2
\biggr)^2 \nonumber\\[-9pt]\\[-9pt]
&&\qquad\quad{} + \mathcal O_p(n \|f_{j,\eta}\|^3).\nonumber
\end{eqnarray}

We first discuss the size of the terms in
(\ref{eq:line3}) and (\ref{eq:line4}). The first term in (\ref
{eq:line3}) can be bounded using the
Cauchy--Schwarz inequality
\[
\sum_{i=1}^n \biggl(\int_{0}^1 (f_{j, \eta}\star g)(x)f_{j,\eta
}(x-\alpha_i)
\,dx \biggr)^2 \leq n \| f_{j, \eta}\star g \|^2 \|f_{j, \eta}\|^2 =
\mathcal O_p(n
\|f_{j,\eta}\|^4).
\]
But observe that
$
\sum_{i=1}^n \mathbb{E}_{W_i} (\int_{0}^1 (f_{j, \eta}\star
g)(x) \,dW_i(x)
)^2= n \|f_{j,\eta} \star g\|^2
$
which does not converge to $0$. Then the Jensen inequality implies
\begin{eqnarray*}
&&\sum_{i=1}^n \mathbb{E}_{W_i} \biggl(\int_{\mathbb{R}} g(\tau
_i) \biggl[\int_{0}^1
f_{j,\eta}(x-\tau_i) \,dW_i(x) \biggr]^2 \,d\tau_i \biggr)^2 \\
&&\qquad\leq\sum_{i=1}^n
\int_{\mathbb{R}} g(\tau_i)\mathbb{E}_{W_i} \biggl[\int_{0}^1
f_{j,\eta}(x-\tau_i)
\,dW_i(x) \biggr]^4 \,d\tau_i
\\
&&\qquad= \mathcal O_p(n \|f_{j,\eta}\|^4).
\end{eqnarray*}
Let us now study the terms derived from double products in
(\ref{eq:line3}) and (\ref{eq:line4}), use first that $2
|ab| \leq
(a^2+b^2)$ to get
$
\sum_{i=1}^n \mathbb{E}_{\alpha_i,W_i} |\int_{0}^1 (f_{j, \eta
}\star
g)(x) f_{j,\eta}(x-\alpha_i) \,dx |
|\int_{\mathbb{R}}
g(\tau_i) [\int_{0}^1 f_{j,\eta}(x-\tau_i) \,dW_i(x) ]^2
\,d\tau_i |
= \mathcal O_p(n \|f_{j,\eta}\|^4)$.

The Cauchy--Schwarz and Jensen inequalities imply
\begin{eqnarray*}
&&\sum_{i=1}^n \mathbb{E}_{W_i} \biggl|\int_{0}^1 (f_{j, \eta}\star
g)(x) \,dW_i(x)
\biggr| \biggl|\int_{\mathbb{R}} g(\tau_i) \biggl[\int_{0}^1
f_{j,\eta}(x-\tau_i)
\,dW_i(x) \biggr]^2 \,d\tau_i \biggr| \\
&&\qquad= \mathcal O_p(n \|f_{j,\eta}\|^3).
\end{eqnarray*}
At last, the Cauchy--Schwarz and Jensen inequalities on the remaining
double-product term imply also
\begin{eqnarray*}
&&{\mathbb E}_{W,\alpha} \Biggl| \sum_{i=1}^n \int_{0}^1 (f_{j, \eta
}\star
g)(x)f_{j,\eta}(x-\alpha_i) \,dx
\int_{0}^1 (f_{j, \eta}\star g)(x) \,dW_i(x) \Biggr| \\
&&\qquad= \mathcal O_p(n \|f_{j,\eta}\|^3).
\end{eqnarray*}
All the above bounds enables us to write
\begin{eqnarray*}
L_1:\!&=&\log\Lambda_n(f_{j, \eta},f_0) \\
&=& \sum_{i=1}^n \int_{0}^1 (f_{j,
\eta}\star g)(x)f_{j,\eta}(x-\alpha_i) \,dx + \int_{0}^1 (f_{j, \eta
}\star g)(x)
\,dW_i(x)\\
&&{} + \frac{1}{2} \sum_{i=1}^n \int_{\mathbb{R}} g(\tau_i)
\biggl[\int_{0}^1
f_{j,\eta}(x-\tau_i) \,dW_i(x) \biggr]^2 - \frac{n}{2} \|f_{j,\eta}\|
^2 \\
&&{} - \frac{1}{2}\sum_{i=1}^n \biggl(\int_{0}^1 (f_{j, \eta}\star g)(x)
\,dW_i(x) \biggr)^2+ \mathcal O_p(n \|f_{j,\eta}\|^3 ).
\end{eqnarray*}
In a similar way, we can also write
\begin{eqnarray*}
L_2:\!&=&\log\Lambda_n(f_{j, \eta^i},f_0)\\
&=& \sum_{i=1}^n \int_{0}^1 (f_{j,
\eta^i}\star g)(x)f_{j,\eta}(x-\alpha_i) \,dx  \\
&&{}+\int_{0}^1 (f_{j, \eta^i}\star
g)(x) \,dW_i(x)\\
&&{} + \frac{1}{2} \sum_{i=1}^n \int_{\mathbb{R}} g(\tau_i)
\biggl[\int_{0}^1
f_{j,\eta^i}(x-\tau_i) \,dW_i(x) \biggr]^2 - \frac{n}{2} \|f_{j,\eta
^i}\|^2 \\
&&{}-\frac{1}{2}\sum_{i=1}^n \biggl(\int_{0}^1 (f_{j, \eta^i}\star
g)(x) \,dW_i(x) \biggr)^2+ \mathcal O_p(n \|f_{j,\eta}\|^3 ).
\end{eqnarray*}
For sake of simplicity, let us write $h=f_{j,\eta^i} - f_{j,\eta} = 2
\eta_i
\psi_{j,i}$. The difference $L = L_2 - L_1$ can thus be decomposed as
%
%
\begin{eqnarray}\qquad
\label{eq:line1.0}
L & = &\sum_{i=1}^n \int_{0}^1 (h \star
g)(x)[f_{j,\eta}(x-\alpha_i)-f_{j,\eta} \star g (x)] \,dx \\
\label{eq:line1.1}
&&{}+ \sum_{i=1}^n \int_{0}^1 (h \star g)(x)(f_{j,\eta} \star g) (x) \,dx +
\int_{0}^1
(h\star g)(x) \,dW_i(x)\\
\label{eq:line1.2}
&&{} + \frac{1}{2} \sum_{i=1}^n \int_{\mathbb{R}} g(\tau_i)
\biggl[\int_{0}^1
f_{j,\eta^i}(x-\tau_i) \,dW_i(x) \biggr]^2 - \frac{n}{2} \|f_{j,\eta
^i}\|^2\\
\label{eq:line1.3}
&&{} - \frac{1}{2} \sum_{i=1}^n \int_{\mathbb{R}} g(\tau_i)
\biggl[\int_{0}^1
f_{j,\eta}(x-\tau_i) \,dW_i(x) \biggr]^2 + \frac{n}{2} \|f_{j,\eta}\|
^2\\
\label{eq:line1.4}
&&{} -  \frac{1}{2} \Biggl[\sum_{i=1}^n \biggl(\int_{0}^1 (f_{j, \eta
^i}\star g)(x)
\,dW_i(x) \biggr)^2 - n \|f_{j,\eta^i} \star g \|^2 \Biggr]\\
\label{eq:line1.4bis}
&&{} - \frac{n}{2} \|f_{j,\eta^i} \star g \|^2
\\
%
%
\label{eq:line1.5}
&&{} + \frac{1}{2} \Biggl[\sum_{i=1}^n
\biggl(\int_{0}^1 (f_{j, \eta}\star g)(x) \,dW_i(x) \biggr)^2 - n \|
f_{j,\eta}
\star g \|^2 \Biggr]\\
\label{eq:line1.5bis}
&&{} + \frac{n}{2}
\|f_{j,\eta} \star g \|^2\\
\label{eq:line1.6}
&&{} + \mathcal O_p(n \|f_{j,\eta}\|^3 ).
\end{eqnarray}

\textit{Bound for} (\ref{eq:line1.0}): we use the classical
Bennett inequality (see, e.g., \cite{massart}) for a sum of
independent and bounded variables. Define
$S = \sum_{i=1}^n \int_{0}^1 (h \star g)(x)[f_{j,\eta}(x-\alpha
_i)-f_{j,\eta}
\star g (x)] \,dx$.
From the Cauchy--Schwarz
inequality, the random variables $ \int_{0}^1 (h \star
g)(x) f_{j,\eta}(x-\alpha_i) \,dx$ are bounded by a constant $b$ such that
$b = \|h \star g\| \|f_{j,\eta}\|$.
Let $v$ and $c$ to be defined as
\[
v = \sum_{i=1}^n \mathbb{E} \biggl[\int_{0}^1 (h \star g) (x)
f_{j,\eta}(x-\alpha_i) \,dx \biggr]^2 \quad\mbox{and}\quad c=b/3.
\]
From the Cauchy--Schwarz inequality, we have that $v \leq n \|f_{j,\eta
}\|^2 \|h \star g\|^2$
and as $h = f_{j,\eta^i}-f_{j,\eta}$, by using our algebraic settings
in Section \ref{sec:alsetting-reg}, we observe that $v \to0$. Bennett's
inequality therefore implies that for any $\kappa>0$
\[
{\mathbb P}(|S| \geq\kappa) \leq2 e^{ {-\kappa
^2}/({2 (n
\|f_{j,\eta}\|^2 \|h \star g\|^2+ \kappa\|h \star g\| \|f_{j,\eta}\|
/3 )})}.
\]
From our algebraic settings in Section \ref{sec:alsetting-reg}, one
has thus that as $n \to\infty$, the ${\mathbb P}(|S| \geq\kappa) $
converges to $0$.

\textit{Bound for} (\ref{eq:line1.2}), (\ref{eq:line1.3}), (\ref
{eq:line1.4}), (\ref{eq:line1.5}): applying Lemma \ref{lem:chi2}
(proved below) to the
chi-square statistics in the expressions (\ref{eq:line1.2}), (\ref{eq:line1.3})
yields that for any $\kappa> 0$
${\mathbb P} ( | \frac{1}{2} \sum_{i=1}^n \int_{\mathbb{R}}
g(\tau_i)
[\int_{0}^1 f_{j,\eta^i}(x-\tau_i) \,dW_i(x) ]^2 - \frac{n}{2}
\|f_{j,\eta^i}\|^2 | \geq\kappa) \leq2
e^{ {-
\kappa^2}/({n \|f_{j,\eta^i}\|^4+ 2 \kappa\|f_{j,\eta^i}\|^2})}$
and
${\mathbb P} ( | \frac{1}{2} \sum_{i=1}^n \int_{\mathbb{R}}
g(\tau_i)
[\int_{0}^1 f_{j,\eta}(x-\tau_i) \,dW_i(x) ]^2 - \frac{n}{2}\times\break
\|f_{j,\eta}\|^2 | \geq\kappa) \leq2 e^{ {-
\kappa^2}/({n \|f_{j,\eta}\|^4+ 2 \kappa\|f_{j,\eta}\|^2})}$.
Similarly, we obtain for the chi-square statistics in (\ref
{eq:line1.4}), (\ref{eq:line1.5}) that for any $\kappa> 0$
\begin{eqnarray*}
&&{\mathbb P} \Biggl( \Biggl| \frac{1}{2} \Biggl[\sum_{i=1}^n
\biggl(\int
_{0}^1 (f_{j, \eta}\star g)(x)
\,dW_i(x) \biggr)^2 - n
\|f_{j,\eta} \star g \|^2 \Biggr] \Biggr| \geq\kappa\Biggr) \\
&&\qquad\leq2 e^{ {- \kappa^2}/({n \|f_{j,\eta}
\star g \|^4+
2 \kappa\|f_{j,\eta} \star g\|^2})}
\end{eqnarray*}
and
\begin{eqnarray*}
&&{\mathbb P} \Biggl( \Biggl| \frac{1}{2} \Biggl[\sum_{i=1}^n
\biggl(\int
_{0}^1 (f_{j, \eta^i}\star g)(x)
\,dW_i(x) \biggr)^2 - n \|f_{j,\eta^i} \star g \|^2 \Biggr] \Biggr|
\geq\kappa\Biggr) \\
&&\qquad\leq2 e^{ {- \kappa^2}/({n \|f_{j,\eta
^i} \star g
\|^4+ 2 \kappa\|f_{j,\eta^i} \star g\|^2})}.
\end{eqnarray*}
It follows from the algebraic setting in Section \ref
{sec:alsetting-reg} that $n\|f_{j,\eta}\|^4 \to0 $ and $\|f_{j,\eta
}\|^2 \to0$, as well as $n \|f_{j,\eta} \star g\|^4 \to0 $ and $\|
f_{j,\eta} \star g\|^2 \to0$ and the above probabilities converge to
zero as $n \to\infty$.

\textit{Bound for} (\ref{eq:line1.1}), (\ref{eq:line1.4bis}), (\ref
{eq:line1.5bis}): using the first term of (\ref{eq:line1.1}), simple
computation shows that
yields
$\sum_{i=1}^n \int_{0}^{1}((f_{j,\eta^i}-f_{j,\eta}) \star g)(x)
(f_{j,\eta}\star g)(x) \,dx - \frac{n}{2} \|f_{j,\eta^i} \star g\|^2 +
\frac{n}{2} \|f_{j,\eta} \star g\|^2 =- \frac{n}{2} \|h \star g\|^2$
and we obtain from our algebraic settings that this term converges to
$0$ since
$n \|h \star g\|^2 \to0$. Moreover, the second term of (\ref
{eq:line1.1}) is
the sum of $n$ i.i.d. centered normal variables and the
Cirelson--Ibragimov--Sudakov inequality \cite{CIS} ensures that
\[
{\mathbb P} \Biggl( \Biggl| \sum_{i=1}^n \int_{0}^1 (h \star g)(x) \,dW_i(x)
\Biggr| \geq
\kappa\Biggr) \leq2 e^{ {{- \kappa^2}/({2 n \|
h\star
g\|^2})}},
\]
and thus the above probability goes to zero.

\textit{Bound for} (\ref{eq:line1.6}): from our algebraic
settings in Section \ref{sec:alsetting-reg}, it follows immediately
that $n \|f_{j,\eta}\|^3 \to0$.

Hence, by combining all the above bounds, it follows that we have shown
that $L_2-L_1$ is the sum of various terms which all converge to zero in
probability or that are larger than some negative constant with
probability tending to one as $n \to+ \infty$, which completes the
proof of Lemma \ref{lem:min-proba}.
\end{pf}
\begin{lemma}\label{lem:chi2}
Let $g$ be a density function on $\mathbb{R}$, and $(W_i)_{i \in\{1,
\ldots, n\}}$ be independent standard Brownian motions on $[0,1]$. Then,
for any $f \in L^{2}([0,1])$ and $\alpha>0$,
\begin{eqnarray*}
&&P \Biggl(\frac{1}{2}\sum_{i=1}^n \int_{\mathbb{R}} g(\tau_i)
\biggl[\int_{0}^1
f(t-\tau_i) \,dW_i(t) \biggr]^2 \,d\tau_i - \frac{n}{2} \|f\|^2 \geq
\alpha\Biggr)
\\
&&\qquad
\leq e^{ {- \alpha^2}/({n \|f\|^4+2 \alpha\|f\|^2})}.
\end{eqnarray*}
\end{lemma}
\begin{pf}
Consider $\zeta_n= \frac{1}{2} \sum_{i=1}^n \int_{\mathbb{R}} g(\tau_i)
[\int_{0}^1 f(t-\tau_i) \,dW_i(t) ]^2 \,d\tau_i - \frac{n}{2} \|f\|^2$.  We
use a Laplace transform technique to bound ${\mathbb P}(\zeta_n \geq
\alpha)$. For any $ \frac{1}{ \|f\|^2}>t>0$, we have by Markov's
inequality
\[
{\mathbb P}(\zeta_n \geq\alpha) \leq e^{ - \alpha
t -
{n}/{2} \|f\|^2 t} \prod_{i=1}^n
\mathbb{E} \bigl[e^{ t/2 \int_{\mathbb{R}}
g(\tau_i) (\int_{0}^1 f(t-\tau_i) \,dW_i(t) )^2 \,d\tau_i} \bigr].
\]
We apply now Jensen's inequality for the exponential function and the measure
$g(\tau) \,d\tau$ to obtain
\[
{\mathbb P}(\zeta_n \geq\alpha) \leq e^{ - \alpha t -
{n}/{2} \|f\|^2 t} \prod_{i=1}^n
\int_{\mathbb{R}} g(\tau_i) \mathbb{E} \bigl[ e^{ t/2
(\int_{0}^1 f(t-\tau_i) \,dW_i(t) )^2} \bigr] \,d\tau_i.
\]
Remark that $ (\int_{0}^1 f(t-\tau_i) \,dW_i(t) )^2$ follows a
chi-square distribution whose Laplace transform does not depend on
$\tau_{i}$ and thus
\[
{\mathbb P}(\zeta_n \geq\alpha) \leq e^{ - \alpha t -
{n}/{2} \|f\|^2 t - {n}/{2} \log
(1 - t \|f\|^2 )}.
\]
Let $ \tilde{\alpha} = \frac{\alpha}{n/2 \|f\|^2}$
and minimizing
now the last bound with respect to $t$ yield the optimal choice $
t^{\star} = \frac{\tilde{\alpha}}{1+\tilde{\alpha}}$. With this
choice, we obtain
${\mathbb P}(\zeta_n \geq\alpha) \leq\break\exp( \frac{n}{2}
[\log(1+\tilde{\alpha}) - \tilde{\alpha} ] )$.
Now use the classical bound
$\log(1+u) - u \leq\frac{-u^2}{2(1+u)}$,
valid for all $u \geq0$, to get
${\mathbb P}(\zeta_n \geq\alpha) \leq\exp( \frac{n}{2} \times
\frac{\tilde{\alpha}^2}{2(1+\tilde{\alpha})} ) = \exp(
\frac{-
\alpha^2}{n \|f\|^4+ 2 \alpha\|f\|^2})$,
which completes the proof of the lemma.
\end{pf}

\subsubsection{A lower bound for the minimax risk}

By Lemmas \ref{lem:assouad-bis} and \ref{lem:min-proba}, it
follows that there exists a constant $C_1$ such
that for all sufficiently large $n$,
\[
\inf_{\hat{f}_n} \sup_{f \in B^s_{p,q}(A)} \mathbb{E} \| \hat
{f}_n - f
\|^2 \geq\inf_{\hat{f}_n} \max_{\eta\in\{\pm1\}^{2^j}}
{\mathbb E}_{f_{j,
\eta}} \|\hat{f}_{n} - f_{j,\eta} \|^{2} \geq C_1 n^{-
({2s})/({2s+2\nu+1})},
\]
which completes the proof of Theorem \ref{theo:lowerbound}.

\subsection{\texorpdfstring{Proof of Theorem
\protect\ref{th:estimshift}}{Proof of Theorem 4}}
For $\bolds\tau= (\tau_{1},\ldots,\tau_{n}) \in{\mathcal T}^{n} $
define the criterion
$M(\bolds\tau) = \frac{1}{n} \sum_{m=1}^{n} \sum_{|\ell| \leq\ell_{0}
} | \theta_{\ell} e^{2 i \ell\pi(\tau_{m}- \tau_{m}^{\ast})} -
\frac{1}{n} \sum_{q=1}^{n} \theta_{\ell} e^{2 i \ell\pi(\tau
_{q}- \tau_{q}^{\ast})} |^{2}$.
Then let us first prove the following lemma.
\begin{lemma} \label{lem:unicity}
Suppose that Assumption \ref{ass:indent} hold. Then, the function
$\bolds\tau\mapsto M(\bolds\tau)$ has a unique minimum on
$\overline{{\mathcal T}
}_{n}$ at $\bolds\tau= \tilde{\bolds\tau}$ such that $M(\tilde
{\bolds\tau}) = 0$
given by
$\tilde{\bolds\tau} = (\tau_{1}^{\ast}-\overline{\tau}_{n},\tau
_{2}^{\ast
}-\overline{\tau}_{n},\ldots,\tau_{n}^{\ast}-\overline{\tau}_{n})$,
where $\overline{\tau}_{n} = \frac{1}{n} \sum_{m=1}^{n} \tau_{m}^{\ast}$.
\end{lemma}
\begin{pf}
By definition of $M(\bolds\tau)$ it follows immediately that
$M(\tilde{\bolds\tau}) = 0$ and thus $\tilde{\bolds\tau}$ is a minimum
since $M(\bolds\tau) \geq0$ for all $\bolds\tau\in{\mathcal T}^{n}$.
Now suppose that there exists $\bolds\tau\in\overline{{\mathcal
T}}_{n}$ such that $M(\bolds\tau)= 0$. This implies that for all
$m=1,\ldots,n$ and all $-\ell_{0} \leq\ell\leq\ell_{0} |
\theta_{\ell}|^{2} | e^{2 i \ell\pi(\tau_{m}- \tau_{m}^{\ast})} -
\frac{1}{n} \sum_{q=1}^{n} e^{2 i \ell\pi(\tau_{q}- \tau_{q}^{\ast})}
|^{2} = 0$.  Since by assumption $\theta_{1}^{\ast} \neq0$, it follows
that for $\ell= 1$, $| e^{2 i \pi(\tau_{m}- \tau_{m}^{\ast})} - \frac
{1}{n} \sum_{q=1}^{n} e^{2 i \pi(\tau_{q}- \tau_{q}^{\ast})} |^{2} = 0$
for all $m=1,\ldots,n$, which implies that $ e^{2 i \pi(\tau_{m}-
\tau_{m}^{\ast})} = e^{2 i \pi(\tau_{q}- \tau_{q}^{\ast})} \mbox{ for
all } m,q=1,\ldots,n$, since $\frac{1}{n} \sum_{q=1}^{n} e^{2 i
\pi(\tau_{q}- \tau_{q}^{\ast})}$ does not depend on $m$. This implies
that $\tau_{m}- \tau_{m}^{\ast} = \tau_{0} \operatorname{mod} 1$
for $m=2,\ldots,n$, where $\tau_{0} = \tau_{1}- \tau_{1}^{\ast}$. By
assumption $\tau_{1},\tau_{1}^{\ast}$ belong to $ {\mathcal T}=
[-\frac{1}{4},\frac{1}{4}]$ and thus $|\tau_{0}| \leq\frac{1}{2}$.
Hence, $\tau_{m} = \tau_{m}^{\ast} + \tau_{0}$ for $m=1,\ldots,n$.
Since $\sum_{m=1}^{n} \tau_{m} = 0$ this implies that $\tau_{0} = -
\frac{1}{n} \sum_{m=1}^{n} \tau_{m}^{\ast}$ and thus $\tau_{m} =
\tilde{\tau}_{m}$ for $m=1,\ldots,n$ which completes the proof.
\end{pf}

Let $F \dvtx{\mathbb{R}}^{n-1} \to{\mathbb{R}}^{n}$ given by
$F(\tau_{2},\ldots,\tau_{n}) = (-\sum_{m=2}^{n} \tau_{m}, \tau
_{2},\ldots,\tau_{n})^{t}$,
and let $\tilde{M} \dvtx{\mathcal T}^{n-1} \to{\mathbb{R}}^{+}$ defined
by $\tilde
{M}(\tau_{2},\ldots,\tau_{n}) = M(F(\tau_{2},\ldots,\tau_{n}))$.
\begin{lemma} \label{lem:hess}
Let $\nabla^{2} \tilde{M}( \tilde{\bolds\tau}_{-1})$ denotes the Hessian
of $\tilde{M}$ at $ \tilde{\bolds\tau}_{-1} = (\tilde{\tau
}_{2},\ldots
,\tilde{\tau}_{n})$, then
$\nabla^{2} \tilde{M}( \tilde{\bolds\tau}_{-1}) = ( \frac{2}{n}
\sum_{|\ell| \leq\ell_{0} } |2 \pi\ell|^{2} |\theta_{\ell}
|^{2} ) ( I_{n-1} + \mathbh{1}_{n-1}^{t} \mathbh{1}_{n-1}
)$,
where $I_{n-1}$ is the $(n-1) \times(n-1)$ identity matrix and
$\mathbh{1}
_{n-1} = (1,\ldots,1)^{t}$ is the vector of ${\mathbb{R}}^{n-1}$ with all
entries equal to one.
Moreover,
$\lambda_{\min}(\nabla^{2} \tilde{M}( \tilde{\bolds\tau}_{-1})) =\break
{\frac
{2}{n} \sum_{|\ell| \leq\ell_{0} } }|2 \pi\ell|^{2} |\theta
_{\ell} |^{2}$,
where $\lambda_{\min}(A)$ denotes the smallest eigenvalue of a
symmetric matrix $A$.
\end{lemma}
\begin{pf}
First remark that for $\bolds\tau_{-1} = (\tau_{2},\ldots,\tau_{n})
\in{\mathbb{R}}^{n-1}$ then $\nabla^{2} \tilde{M}(\bolds\tau_{-1}) =
\nabla F ^{t} \nabla^{2} M(F(\bolds\tau_{-1})) \nabla F $ where
$\nabla^{2} M(F(\bolds\tau_{-1}))$ denotes the Hessian of $M$ at $
F(\bolds\tau_{-1}) $ and $\nabla F$ is the gradient of $F$ [$n
\times(n-1)$ matrix not depending on $\tau$].
Now, since for any $\bolds\tau\in{\mathcal T}^{n}$ $
M(\bolds\tau) = {\sum_{|\ell| \leq\ell_{0} } }| \theta_{\ell}|^{2}
(1- | \frac{1}{n} \sum_{q=1}^{n} e^{2 i \ell\pi(\tau
_{q}- \tau_{q}^{\ast})} |^{2} )
$
it follows that for $m=2,\ldots,n$
\[
\frac{\partial}{\partial\tau_{m}} M(\bolds\tau) = {-\frac{2}{n^2}
\sum
_{|\ell| \leq\ell_{0} }} | \theta_{\ell}|^{2} \Re\Biggl[ 2 i \pi
\ell e^{2 i \ell\pi(\tau_{m}- \tau_{m}^{\ast})} \Biggl( \overline
{ \sum_{q=1}^{n} e^{2 i \ell\pi(\tau_{q}- \tau_{q}^{\ast})} }
\Biggr) \Biggr],
\]
where $\Re[ z ]$ denotes the real part of a complex number. Hence, for
$m_{1} \neq m_{2}$
\[
\frac{\partial^2}{\partial\tau_{m_2}\, \partial\tau_{m_1}} M(\bolds
\tau
) = {-\frac{2}{n^2} \sum_{|\ell| \leq\ell_{0} } }|2 \pi\ell|^{2} |
\theta_{\ell}|^{2} \Re\bigl[ e^{2 i \ell\pi(\tau_{m_1}- \tau
_{m_1}^{\ast} - \tau_{m_2} + \tau_{m_2}^{\ast} ) } \bigr]
\]
and for $m_1 = m_2$
\begin{eqnarray*}
&&\frac{\partial^2}{\partial\tau_{m_1}\, \partial\tau_{m_1}} M(\bolds
\tau)\\
&&\qquad = -\frac{2}{n^2} \sum_{|\ell| \leq\ell_{0} } |2 \pi\ell|^{2} |
\theta_{\ell}|^{2} \Re\Biggl[1 - e^{2 i \ell\pi(\tau_{m_1}- \tau
_{m_1}^{\ast})} \Biggl( \overline{ \sum_{q=1}^{n} e^{2 i \ell\pi
(\tau_{q}- \tau_{q}^{\ast})} } \Biggr) \Biggr].
\end{eqnarray*}
Then, remark that $F(\tilde{\bolds\tau}_{-1}) = \tilde{\bolds\tau
}$. Hence,
by taking $\tau_{m} = \tilde{\tau}_{m}$ for $m=2,\ldots,n$ in the
above formulas, it follows that
%
%
\begin{equation} \label{eq:gradM}\qquad
\nabla^{2} M( \tilde{\bolds\tau} ) = \nabla^{2} M(F(\tilde{\bolds
\tau
}_{-1})) = {\frac{2}{n} \sum_{|\ell| \leq\ell_{0} }} |2 \pi\ell
|^{2} |\theta_{\ell} |^{2} \biggl( I_{n} - \frac{1}{n} \mathbh
{1}_{n}^{t} \mathbh{1}
_{n} \biggr),
\end{equation}
where $I_{n}$ is the $n \times n$ identity matrix and $\mathbh{1}_{n} =
(1,\ldots,1)^{t}$ is the vector of ${\mathbb{R}}^{n}$ with all
entries equal
to one. Hence, the result follows from (\ref{eq:gradM}) and the
equality $\nabla^{2} \tilde{M}(\tilde{\bolds\tau}_{-1}) = \nabla F ^{t}
\nabla^{2} M(F(\tilde{\bolds\tau}_{-1})) \nabla F $, and the fact that
the eigenvalues of the matrix $A = I_{n-1} + \mathbh{1}_{n-1}^{t}
\mathbh{1}_{n-1}$ are
$n$ (of multiplicity 1) and $1$ (of multiplicity $n-2$).
\end{pf}
\begin{lemma} \label{lem:quadra}
Suppose that Assumption \ref{ass:indent} holds. Then, there exists a
constant $\kappa(f) >0$ (depending on the shape function $f$) such
that for all $\bolds\tau\in\overline{{\mathcal T}}_{n}
M(\bolds\tau) - M(\tilde{\bolds\tau}) \geq\kappa(f) ( {\sum
_{|\ell|
\leq\ell_{0} } }| \theta_{\ell}|^{2} ) ( \frac{1}{n}
\sum_{m=2}^{n} (\tau_{m} -\tilde{\tau}_{m} )^{2} )$.
\end{lemma}
\begin{pf}
First, remark that for any $\bolds\tau\in\overline{{\mathcal T}}_{n}$
then $\tilde{M}(\bolds\tau_{-1}) = M(F(\bolds\tau))$ where
$\bolds\tau_{-1} = (\tau_{2},\ldots,\tau_{n})$. Since
$\tilde{\bolds\tau}$ is a minimum of $\bolds\tau\mapsto M(\bolds\tau)$,
a second order Taylor expansion implies that for all $\bolds\tau_{-1}$
in neighborhood ${\mathcal V}\subset{\mathcal T}^{n-1}$ of $\tilde
{\bolds\tau}_{-1}$
\begin{eqnarray*}
M(\bolds\tau) - M(\tilde{\bolds\tau}) & = & \tilde{M}(\bolds\tau
_{-1}) - \tilde
{M}(\tilde{\bolds\tau}_{-1}) \\
& = & (\bolds\tau_{-1}-\tilde{\bolds\tau}_{-1})^{t} \nabla^{2}
\tilde
{M}(\tilde{\bolds\tau}_{-1}) (\bolds\tau_{-1}-\tilde{\bolds\tau
}_{-1}) + o( \|
\bolds\tau_{-1}-\tilde{\bolds\tau}_{-1} \|^{2}).
\end{eqnarray*}
Using Lemma \ref{lem:hess} and the above equation, it follows that
there exists a universal constant $0 < c_{1} < 1$ and an open
neighborhood $\tilde{{\mathcal V}} \subset{\mathcal V}$ of $\tilde
{\bolds\tau}$ such
that for all $\bolds\tau\in\tilde{{\mathcal V}} $
\[
M(\bolds\tau) - M(\tilde{\bolds\tau}) \geq2 c_{1} \biggl( {\sum
_{|\ell| \leq
\ell_{0} }} |2 \pi\ell|^{2} |\theta_{\ell} |^{2} \biggr)
\Biggl(\frac{1}{n} \sum_{m=2}^{n} (\tau_{m} -\tilde{\tau}_{m} )^{2}
\Biggr).
\]
Now remark that under Assumption \ref{ass:indent}, $M(\bolds\tau) >
M(\tilde{\bolds\tau}) = 0$ for all $\bolds\tau\in\overline
{{\mathcal T}}_{n}
\setminus\tilde{{\mathcal V}}$ by Lemma \ref{lem:unicity}. Since
$M(\bolds\tau) = {\sum_{|\ell| \leq\ell_{0} }} | \theta_{\ell}|^{2}
(1- | \frac{1}{n} \sum_{q=1}^{n} e^{2 i \ell\pi(\tau
_{q}- \tau_{q}^{\ast})} |^{2} )$,
the compactness of $\overline{{\mathcal T}}_{n}$ and the continuity of
$\tau
\mapsto M(\bolds\tau)$ implies that there exists a constant $0 <
c_{2}(f) <
1$ (depending on $\tilde{{\mathcal V}}$ and thus on $f$) such that for all
$\bolds\tau\in\overline{{\mathcal T}}_{n} \setminus\tilde
{{\mathcal V}}$,
$M(\bolds\tau) \geq{\sum_{|\ell| \leq\ell_{0} } }| \theta_{\ell}|^{2}
(1-c_{2}(f) )$.
Moreover, since ${\mathcal T}$ is a compact set it follows that there
exists a
universal constant $c_{3} > 0$ such that $(\tau_{m} -\tilde{\tau
}_{m} )^{2} \leq c_{3}$ for all $m=2,\ldots,n$, which implies that for
all $\bolds\tau\in\overline{{\mathcal T}}_{n}
\frac{1}{n}\sum_{m=2}^{n} (\tau_{m} -\tilde{\tau}_{m} )^{2} \leq c_{3}$.
Therefore,
\[
M(\bolds\tau) - M(\tilde{\bolds\tau}) \geq\biggl( c_{3}^{-1}
\bigl(1-
c_{2}(f) \bigr) \sum_{|\ell| \leq\ell_{0} } | \theta_{\ell}|^{2}
\biggr) \Biggl( \frac{1}{n} \sum_{m=2}^{n} (\tau_{m} -\tilde{\tau
}_{m} )^{2} \Biggr)
\]
for all $\bolds\tau\in\overline{{\mathcal T}}_{n} \setminus\tilde
{{\mathcal V}}$. Then
the result follows by taking $\kappa(f) = \min(2 c_{1}, c_{3}^{-1}
(1- c_{2}(f) ))$ and the fact that ${\sum_{|\ell| \leq
\ell_{0} }} |2 \pi\ell|^{2} |\theta_{\ell} |^{2} \geq{\sum
_{|\ell| \leq\ell_{0} } }| \theta_{\ell}|^{2}$.
\end{pf}

Now recall that
$\hat{\bolds\tau} = (\hat{\tau}_{1},\ldots,\hat{\tau}_{n}) =
\arg\min
_{ \bolds\tau\in\overline{{\mathcal T}}_{n} } M_{n}(\bolds\tau)$.
Since $\hat{\bolds\tau}$ is a minimum of $\bolds\tau\mapsto
M_{n}(\bolds\tau)$ and
$\tilde{\bolds\tau}$ is a minimum of $\bolds\tau\mapsto M(\bolds
\tau)$, it follows that
$M(\hat{\bolds\tau}) - M(\tilde{\bolds\tau}) \leq{2 \sup_{\bolds
\tau\in{\mathcal T}^{n}}}
|M_{n}(\bolds\tau) - M(\bolds\tau)|$.
Therefore, Lemma \ref{lem:quadra} implies that
%
%
\begin{equation} \label{eq:tau}\qquad
\frac{1}{n} \sum_{m=2}^{n} (\hat{\tau}_{m} -\tilde{\tau}_{m}
)^{2} \leq2 {\biggl( \kappa(f) \biggl( {\sum_{|\ell| \leq\ell_{0} }}
| \theta_{\ell}|^{2} \biggr) \biggr)^{-1} \sup_{\bolds\tau\in
{\mathcal T}
^{n}}} |M_{n}(\bolds\tau) - M(\bolds\tau)|.
\end{equation}
\begin{lemma} \label{lem:concentration}
Let $Z = {\sup_{\bolds\tau\in{\mathcal T}^{n}}} |M_{n}(\bolds\tau)
- M(\bolds\tau)|$. Then
for any $t > 0$
\[
{\mathbb P} \biggl( Z \leq\biggl(1+ 2 \biggl({ \sum_{|\ell| \leq\ell
_{0} }}
| \theta_{\ell}|^{2} \biggr)^{1/2} \biggr) \bigl( \sqrt
{C(\varepsilon,n,\ell_{0},t)} + C(\varepsilon,n,\ell_{0},t) \bigr)
\biggr) \geq1 - \exp(-t),
\]
where $C(\varepsilon,n,\ell_{0},t) = \varepsilon^2 (2 \ell_{0} + 1) + 2
\varepsilon^2 \sqrt{\frac{2 \ell_{0} + 1}{n}t} + 2
\frac{\varepsilon^2}{n}t$.
\end{lemma}
\begin{pf}
Remark that $M_{n}(\bolds\tau)$ can be
decomposed as
$M_{n}(\bolds\tau) = M(\bolds\tau) + L(\bolds\tau) + Q(\bolds\tau)$,
where
\begin{eqnarray*}
L(\bolds\tau) & = & 2 \frac{\varepsilon}{n} \sum_{m=1}^{n} \sum
_{|\ell|
\leq\ell_{0} } \Re\Biggl[ \Biggl( \theta_{\ell} e^{2 i \ell\pi
(\tau_{m}- \tau_{m}^{\ast})} - \frac{1}{n} \sum_{q=1}^{n} \theta
_{\ell} e^{2 i \ell\pi(\tau_{q}- \tau_{q}^{\ast})} \Biggr)\\
&&\hspace*{84.95pt}{}\times
\Biggl( \overline{z_{m,\ell} e^{2 i \ell\pi\tau_{m}} - \frac{1}{n} \sum
_{q=1}^{n} z_{q,\ell} e^{2 i \ell\pi\tau_{q}}} \Biggr) \Biggr], \\
Q(\bolds\tau) & = & \frac{\varepsilon^2}{n} \sum_{m=1}^{n} \sum
_{|\ell|
\leq\ell_{0} } \Biggl|z_{m,\ell} e^{2 i \ell\pi\tau_{m}} - \frac
{1}{n} \sum_{q=1}^{n} z_{q,\ell} e^{2 i \ell\pi\tau_{q}} \Biggr|^{2}.
\end{eqnarray*}
By the Cauchy--Schwarz inequality, $| L(\bolds\tau) | \leq2 \sqrt
{M(\bolds\tau)}
\sqrt{Q(\bolds\tau)} $. Since $M(\bolds\tau) \leq\break{\sum_{|\ell|
\leq\ell
_{0} } }| \theta_{\ell}|^{2}$ for all $\bolds\tau\in{\mathcal
T}^{n}$, one has that\vspace*{1pt}
$| L(\bolds\tau) | \leq2 ( {\sum_{|\ell| \leq\ell_{0} }} |
\theta
_{\ell}|^{2} )^{1/2} \sqrt{Q(\bolds\tau)}$.
Therefore,
%
%
\begin{eqnarray} \label{eq:Delta}
&&{\sup_{\bolds\tau\in{\mathcal T}^{n}}} | M_{n}(\bolds\tau) -
M(\bolds\tau) | \nonumber\\[-8pt]\\[-8pt]
&&\qquad\leq\biggl(1+
2 \biggl( \sum_{|\ell| \leq\ell_{0} } | \theta_{\ell}|^{2}
\biggr)^{1/2} \biggr) \biggl( \sqrt{ \sup_{\bolds\tau\in{\mathcal
T}^{n}} Q(\bolds\tau)}
+ \sup_{\bolds\tau\in{\mathcal T}^{n}} Q(\bolds\tau) \biggr).\nonumber
\end{eqnarray}
Thus, it suffices to derive a concentration inequality for $\sup
_{\bolds\tau\in{\mathcal T}^{n}} Q(\bolds\tau)$. For this remark that
$
Q(\bolds\tau) \leq W_{1}
$
for all $\bolds\tau\in{\mathcal T}^{n}$, where $W_{1} = {\sum_{|\ell
| \leq\ell
_{0} } \frac{\varepsilon^2}{n} \sum_{m=1}^{n}} |z_{m,\ell}|^{2} $.
Then using a standard concentration inequality for sum of $\chi^2$
variables (see, e.g., \cite{massart}) one has that for any $t > 0$
${\mathbb P} ( \sup_{\bolds\tau\in{\mathcal T}^{n}} Q(\bolds
\tau) \geq C(\varepsilon
,n,\ell_{0},t) ) \leq{\mathbb P} (W_{1} \geq C(\varepsilon
,n,\ell
_{0},t) ) \leq\exp(-t)$,
where $C(\varepsilon,n,\ell_{0},t) = \varepsilon^2 (2 \ell_{0} + 1) +\break 2
\varepsilon^2 \sqrt{\frac{2 \ell_{0} + 1}{n}t} + 2 \frac{\varepsilon
^2}{n}t$. Therefore, the result follows using inequality (\ref
{eq:Delta}).
\end{pf}
%

From Lemma \ref{lem:concentration} and inequality (\ref{eq:tau}), it
follows that
%
%
\begin{eqnarray}\label{eq:tauproba}
&&{\mathbb P} \Biggl( \frac{1}{n} \sum_{m=2}^{n} (\hat{\tau}_{m}
-\tilde
{\tau}_{m} )^{2}\leq\frac{ 2+ 4 ( \sum_{|\ell| \leq\ell
_{0} } | \theta_{\ell}|^{2} )^{1/2}}{ \kappa(f) ( \sum
_{|\ell| \leq\ell_{0} } | \theta_{\ell}|^{2} ) } \nonumber\\
&&\hspace*{100.5pt}{}\times\bigl(
\sqrt{C(\varepsilon,n,\ell_{0},t)} + C(\varepsilon,n,\ell_{0},t)
\bigr) \Biggr) \\
&&\qquad\geq1 - \exp(-t).\nonumber
\end{eqnarray}
To complete the proof, remark that
$\frac{1}{n} \sum_{m=2}^{n} (\hat{\tau}_{m}-\tau_{m}^{\ast} )^{2}
\leq2 ( \frac{1}{n} \sum_{m=2}^{n} (\hat{\tau}_{m}-\tilde
{\tau}_{m} )^{2}
+ (\frac{1}{n} \sum_{m=1}^{n} \tau_{m}^{\ast}
)^{2} )$.
Since the $\tau_{m}^{\ast}$ are i.i.d. variables with zero mean and
bounded by $1/4$, Bernstein's inequality (see, e.g., \cite{massart})
implies that for any $t > 0$ then
%
%
\begin{equation} \label{eq:meantau}
{\mathbb P} \Biggl( \Biggl|\frac{1}{n} \sum_{m=1}^{n} \tau_{m}^{\ast}
\Biggr| \geq\sqrt{2 \sigma^{2}_{g} \frac{t}{n} } + \frac{t}{12 n}
\Biggr) \leq2 \exp(-t),
\end{equation}
where $ \sigma^{2}_{g} = \int_{{\mathcal T}} \tau^{2} g(\tau) \,d \tau$.
Then Theorem \ref{th:estimshift} follows from inequalities (\ref
{eq:tauproba}) and (\ref{eq:meantau}).

\subsection{\texorpdfstring{Proof of Theorem
\protect\ref{theo:VanTree}}{Proof of Theorem 5}}

To simplify the notation, we write $\tau_{m} = \tau_{m}^{\ast}$ to
denote the true shifts. Part of the proof is inspired by general
results on Van Tree inequalities in \cite{gillev}. First, let us
consider the case where the shifts $ \tau_m, m=1,\ldots,n$, are
fixed parameters to estimate and let $\tau^{n} = (\tau_1,\ldots,\tau
_{n})$. Recall that $X = (c_{m,\ell})_{\ell\in{\mathbb{Z}},
m=1,\ldots, n}$
denote the set of observations taking values in the set ${\mathcal X}=
{\mathbb{C}}
^{\infty\times n}$. Then, the likelihood of the random variable $X$ is
given by
$p(x | \tau^{n}) = C\prod_{m=1}^{n} \prod_{\ell\in{\mathbb{Z}}}
\exp
\{ -\frac{1}{2 \varepsilon^{2}} |c_{m,\ell} - \theta_{\ell}
e^{-i 2 \pi\ell\tau_m} |^{2} \}$.
Therefore, for $m=1,\ldots,n$,
%
%
\begin{equation} \label{eq:score}
{\mathbb E}_{\tau} \biggl( \frac{\partial}{\partial\tau_{m}} \log
p(x |
\tau^{n}) \biggr) = 0,
\end{equation}
where for a function $h(X)$ of the random variable $X$, ${\mathbb
E}_{\tau}
h(X) = \int_{{\mathcal X}} h(x) p(x |\break \tau^{n}) \,dx$.
Then, for $m_{1} \neq m_{2}$ one has that
${\mathbb E}_{\tau} ( \frac{\partial}{\partial\tau_{m_{1}}}
\log p(x
| \tau^{n}) \frac{\partial}{\partial\tau_{m_{2}}} \log p(x | \tau
^{n}) ) = 0$,
and for $m_{1} = m_{2}
{\mathbb E}_{\tau} ( \frac{\partial}{\partial\tau_{m_{1}}}
\log p(x
| \tau^{n}) )^{2} = \varepsilon^{-2} \sum_{\ell\in{\mathbb
{Z}}} (2\pi
\ell)^{2} |\theta_{\ell}|^{2}$.

Now assume that the shifts are i.i.d. random variables with density
$g(\tau)$ satisfying Assumption \ref{ass:gT}. Let $\hat{\tau}^{n}
= \hat{\tau}^{n}(X)$ denote any estimator of the shifts $\tau^{n}$.
Then define the following vectors $U$ and $V = (V_{1},\ldots,V_{n})'$
in ${\mathbb{R}}^{n}$ as
\[
U = \hat{\tau}^{n} - \tau^{n}
\]
and
\[
V_{m} = \frac{\partial
}{\partial\tau_{m}} [ p(x | \tau^{n}) g_{n}(\tau^{n}) ] \frac{1}{
p(x | \tau^{n}) g_{n}(\tau^{n})} \qquad\mbox{for } m=1,\ldots,n,
\]
where $g_{n}(\tau^{n}) = \prod_{m=1}^{n}g(\tau_m)$. First, remark that
\begin{eqnarray*}
{\mathbb E} ( U'V ) & = & \int_{{\mathcal X}} \int
_{{\mathcal T}^{n}} \sum
_{m=1}^{n} (\hat{\tau}^{n}_{m} - \tau_{m}) \,\frac{\partial
}{\partial\tau_{m}} [ p(x | \tau^{n}) g_{n}(\tau^{n}) ] \,d \tau^{n}\,
d x \\
& = & \int_{{\mathcal X}} \sum_{m=1}^{n} \hat{\tau}^{n}_{m} \biggl(
\int
_{{\mathcal T}^{n}} \frac{\partial}{\partial\tau_{m}} [ p(x | \tau^{n})
g_{n}(\tau^{n}) ] \,d \tau^{n} \biggr)\,d x \\
& &{} - \int_{{\mathcal X}} \sum_{m=1}^{n} \biggl( \int_{{\mathcal
T}^{n}} \tau_{m}\,
\frac{\partial}{\partial\tau_{m}} [ p(x | \tau^{n}) g_{n}(\tau
^{n}) ] \,d \tau^{n} \biggr)\,d x.
\end{eqnarray*}
An integration by parts and the fact that $\lim_{\tau\to\tau_{\min
}} g(\tau) = \lim_{\tau\to\tau_{\max}} g(\tau) = 0$ implies that
$\int_{{\mathcal T}^{n}} \frac{\partial}{\partial\tau_{m}} [ p(x |
\tau
^{n}) g_{n}(\tau^{n}) ] \,d \tau^{n} = 0$.
Using again an integration by parts and Assumption \ref{ass:g}, one has that
$\int_{{\mathcal T}^{n}} \tau_{m} \,\frac{\partial}{\partial\tau
_{m}} [ p(x
| \tau^{n}) g_{n}(\tau^{n}) ] \,d \tau^{n} = - \int_{{\mathcal
T}^{n}} p(x |
\tau^{n})\times g_{n}(\tau^{n}) \,d \tau^{n}$.
Therefore,
${\mathbb E} ( U'V ) = \sum_{m=1}^{n} \int_{{\mathcal
T}^{n}} \int_{{\mathcal X}
} p(x | \tau^{n}) g_{n}(\tau^{n}) \,d \tau= n$.\vspace*{1pt}
Now using the Cauchy--Schwarz inequality, it follows that
$n^{2} = ( {\mathbb E} ( U'V ) )^{2} \leq
{\mathbb E} (
U'U )\times {\mathbb E} ( V' V )$.
Then remark that
\[
{\mathbb E} ( U'U ) = {\mathbb E}\sum_{m=1}^{n} (\hat
{\tau}^{n}_{m}-
\tau_{m})^{2} = \int_{X} \int_{{\mathcal T}^{n}} \bigl(\hat{\tau}^{n}_{m}(x)-
\tau_{m}\bigr)^{2} p(x | \tau^{n}) g_{n}(\tau^{n}) \,dx \,d\tau
\]
and
\begin{eqnarray*}
{\mathbb E} ( V'V ) & = & {\mathbb E}\sum_{m=1}^{n}
\biggl( \frac
{\partial}{\partial\tau_{m}} [ \log p(x | \tau^{n}) + \log
g_{n}(\tau^{n}) ] \biggr)^{2} \\
& = & {\mathbb E}\sum_{m=1}^{n} \biggl( \frac{\partial}{\partial
\tau_{m}}
\log p(x | \tau^{n}) \biggr)^{2} + {\mathbb E}\sum_{m=1}^{n} \biggl(
\frac
{\partial}{\partial\tau_{m}} \log g_{n}(\tau^{n}) \biggr)^{2},
\end{eqnarray*}
since by using (\ref{eq:score}) it follows that
$
{\mathbb E} ( \sum_{m=1}^{n} \frac{\partial}{\partial\tau
_{m}} \log
p(x | \tau^{n}) \,\frac{\partial}{\partial\tau_{m}} \log g_{n}(\tau
^{n}) ) = \sum_{m=1}^{n} \int_{{\mathcal T}^{n}} ( \int
_{{\mathcal X}}
( \frac{\partial}{\partial\tau_{m}} \log p(x | \tau^{n})
) p(x | \tau^{n}) \,dx ) ( \frac{\partial
}{\partial\tau_{m}} \log g_{n}(\tau^{n}) ) g_{n}(\tau^{n})
\,d\tau^{n} = 0$.
Hen\-ce,
\begin{eqnarray*}
{\mathbb E} ( V'V ) & = & n \varepsilon^{-2} \sum_{\ell\in
{\mathbb{Z}}}
(2\pi\ell)^{2} |\theta_{\ell}|^{2} + {\mathbb E}\sum_{m=1}^{n}
\biggl(
\frac{\partial}{\partial\tau_{m}} \log g(\tau_{m} ) \biggr)^{2} \\
& = & n \varepsilon^{-2} \sum_{\ell\in{\mathbb{Z}}} (2\pi\ell)^{2}
|\theta
_{\ell}|^{2} + n \int_{{\mathcal T}} \biggl( \frac{\partial
}{\partial\tau
} \log g(\tau) \biggr)^{2} g(\tau) \,d \tau,
\end{eqnarray*}
which completes the proof using that $n^{2} \leq{\mathbb E} ( U'U
) {\mathbb E} ( V'V )$.
\end{appendix}

\section*{Acknowledgments}

We would like to thank Patrick Cattiaux for fruitful discussions. Both
authors are very much indebted to the two anonymous referees and the
Associate Editor for their constructive criticism, comments and remarks
that resulted in a major revision of the original manuscript.

\printaddresses

\end{document}